\documentclass[11pt,leqno]{amsart}
\usepackage{amssymb}
\usepackage{mathrsfs}
\usepackage{color}
\usepackage{soul}
\usepackage{dsfont}

\newcommand\N{{\mathbb N}}
\newcommand\R{{\mathbb R}}

\newcommand\C{{\mathbb C}}

\newcommand\Z{{\mathbb Z}}

\def\AA{{\mathcal A}}
\def\BB{{\mathcal B}}
\def\CC{{\mathcal C}}
\def\DD{{\mathcal D}}

\def\FF{{\mathcal F}}

\def\LL{{\mathcal L}}

\def\OO{{\mathcal O}}

\def\QQ{{\mathcal Q}}

\def\SS{{\mathcal S}}
\def\TT{{\mathcal T}}

\def\VV{{\mathcal V}}

\def\BBB{{\mathscr{B}}}

\def\eps{{\varepsilon}}

\newtheorem{theo}{Theorem}
\newtheorem{prop}[theo]{Proposition}
\newtheorem{lem}[theo]{Lemma}

\newcommand{\beqn}{\begin{equation}}
\newcommand{\eeqn}{\end{equation}}
\newcommand{\bear}{\begin{eqnarray}}
\newcommand{\eear}{\end{eqnarray}}
\newcommand{\bean}{\begin{eqnarray*}}
\newcommand{\eean}{\end{eqnarray*}}
\newcommand{\bal}{\begin{aligned}}
\newcommand{\eal}{\end{aligned}}

\newcommand{\be}{\begin{equation}}
\newcommand{\ee}{\end{equation}}
\newcommand{\ba}{\begin{aligned}}
\newcommand{\ea}{\end{aligned}}

\def\Nt{|\hskip-0.04cm|\hskip-0.04cm|}

\newcommand{\Black}{\color{black}}

\newcommand{\Blue}{\color{black}}


\begin{document}

\address[St\'ephane Mischler]{Universit\'e Paris-Dauphine, Institut Universitaire de France (IUF), PSL Research University,
CNRS, UMR [7534], CEREMADE, 
Place du Mar\'echal de Lattre de Tassigny
75775 Paris Cedex 16, 
France.}
\email{mischler@ceremade.dauphine.fr}

\address[Qilong Weng]{Universit\'e Paris-Dauphine, PSL Research University,
CNRS, UMR [7534], CEREMADE, 
Place du Mar\'echal de Lattre de Tassigny
75775 Paris Cedex 16, 
France.}
\email{weng@ceremade.dauphine.fr}

\title{On a linear runs and tumbles equation}


\author{S. Mischler, Q. Weng}

\begin{abstract}
We consider a  linear runs and tumbles equation in dimension $d \ge 1$ for which 
we establish the existence of a unique positive and normalized steady state as well as its asymptotic stability, 
improving similar results obtained by Calvez et al. \cite{CRS} in dimension $d=1$. Our analysis is based on
the Krein-Rutman theory revisited in \cite{MiSch2016} together with some new moment estimates for proving confinement 
mechanism as well as dispersion, multiplicator and averaging lemma arguments for proving some regularity property of suitable iterated averaging quantities.

\end{abstract}

\maketitle

\begin{center} {\bf Version of \today}
\end{center}

\vspace{0.3cm}

\bigskip
\bigskip


\tableofcontents


\section{Introduction and main result}
\label{sec:model}
\setcounter{equation}{0}
\setcounter{theo}{0}

\subsection{The ``runs and tumbles" equation in chemotaxis}

In the present paper we are interested in a kinetic evolution PDE coming from the modeling
\Blue of cells movement in the presence of a chemotactic chemical substance. 
The so-called {\it run-and-tumble} model introduced by Stroock \cite{Stroock} and Alt \cite{Alt}, and studied further in  \cite{ODA,OthmerStevens,ErbanOthmer}, \Black
\beqn\label{eq:RTeq1}
\partial_t f = \LL f = - v \cdot \nabla_x f + \int_\VV \bigl\{ K' f' - K f \bigr\}  \, dv' 
\eeqn
describes the evolution of the distribution function of a microorganisms density $f=f(t,x,v) \ge 0$ which at time $t \ge0$ and at position $x\in\R^d$ move with the velocity $v \in \VV$. At a microscopic description level, microorganisms move in straight line with their own velocity $v$ which changes accordingly to a jump process of parameter $K = K(x,v,v') \ge 0$. Here and below, we used the shorthands  $f' =f(t,x,v')$   and $K' = K(x,v',v)$. 
For the sake of simplicity, we assume that $\VV \subset \R^d$ is the centered ball with unit volume ($\VV := B(0,V_0) $ with $V_0$ chosen such that  $|\VV| = 1$). 
We complement the evolution PDE \eqref{eq:RTeq1} with an  initial condition 
\beqn\label{eq:InitCond} 
f(0,\cdot) = f_0 \quad \hbox{on}\quad \R^d \times \VV.
\eeqn

\smallskip 
At least formally, for any multiplier $\varphi = \varphi(x,v)$, we have 
$$
{d \over dt} \int f \varphi = \int f  \, v \cdot \nabla_x \varphi + \int f K  \int_\VV \bigl\{ \varphi' - \varphi \bigr\} .
$$
In particular, choosing $\varphi \equiv 1$ in the above identity, we see that the total mass is conserved and we may assume that it is normalized to the unit,
namely 
\beqn\label{eq:Conservation}
\langle\!\langle f(t,\cdot) \rangle\!\rangle = \langle\!\langle f_0 \rangle\!\rangle = 1,
\qquad t \ge 0, 
\eeqn
where for functions $g = g(x,v)$ and $h = h(v)$, we  define
$$
\langle h \rangle = \int_\VV h \, dv, \quad \langle\!\langle g \rangle\!\rangle  = \int_{\R^d} \langle g(x,\cdot ) \rangle \, dx .
$$

\smallskip
The precise form of the turning kernel $K$ depends upon the \Blue possibly time and space dependent \Black concentration $S = S(t,x)$ of a chemical agent: microorganisms have the tendency to move to where the chemical concentration is higher. More specifically,  we assume that the turning kernel is given by 
\beqn\label{eq:RTeq2}
K = K[S](v) := 1 - \chi \, \Phi( \Blue \partial_t S + v \cdot \nabla_x S \Black),  \quad \chi \in (0,1), \,\, \Phi(y) = \hbox{sign} (y),
\eeqn
where the sign function is defined by $\Phi(y) = - 1$ if $y < 0$ and $\Phi(y) = 1$ if $y > 0$. 
In other words, the turning kernel $K$ takes the two values  $1 \pm \chi$ depending on the velocity direction of the microorganism with respect to the gradient of the  chemical concentration.  

\smallskip
When the chemical agent is produced by the microorganisms themselves, it is usually assumed to be given as the solution to the damped Poisson equation 
\beqn\label{eq:RTeq3}
- \Delta S + S = \varrho := \int_\VV f \, dv,
\eeqn
so that the evolution of the  microorganisms density $f$ is given by the coupled system of equations \eqref{eq:RTeq1}-\eqref{eq:RTeq2}-\eqref{eq:RTeq3}. 
We refer the reader interested by the well-posedness issue for related models to the  review paper \cite{MR2809625} and the references quoted therein.
We also refer to \cite{MR3011314} for related modeling considerations. 
Concerning the qualitative behaviour of the solutions it seems that the unique available information is the mass conservation \eqref{eq:Conservation}. 
One of the main difficulty comes from the fact that both equations  \eqref{eq:RTeq1} and \eqref{eq:RTeq3} are invariant by translations so that the expected 
confinement mechanism seems to be hard to prove. 

\smallskip
On the other way round, one can see that for a given solution $(f,S)$ of  \eqref{eq:RTeq1}-\eqref{eq:RTeq2}-\eqref{eq:RTeq3} and for any rotation $R \in SO(d)$ the couple $(f_R,S_R)$ is also a solution of the same equations, where we have set $f_R(x,v) := f(Rx,Rv)$ and $S_R(x) := S(Rx)$. In particular, an invariant by rotations initial datum $f_0$ gives rise to an invariant by rotations solution $(f,S)$.  

\smallskip
More specifically, we may observe that in the case when $S$ does not depend of time and it is radially symmetric and strictly decreasing in the position variable (which is the case if the density $\varrho$ satisfies the same properties thanks to the maximum principle), we have  
$$
- \Phi( \Blue \partial_t S + v \cdot \nabla_x S \Black) = - \Phi( - v \cdot x) =  \hbox{sign} (x \cdot v), 
$$
and thus the associated turning kernel writes 
\beqn\label{eq:RTeq4}
K = K(x,v) := 1 + \chi \zeta, \quad \chi \in (0,1), \,\, \zeta = \zeta(x,v) = \hbox{sign} (x \cdot v).
\eeqn
Such a kernel has been introduced in \cite{CRS} and the  associated (now linear!) evolution equation  \eqref{eq:RTeq1}-\eqref{eq:RTeq4} has then been analyzed in dimension $d=1$: the existence of a unique (positive and normalized) steady state has been established and its asymptotic exponential stability has been proved. 

\subsection{The linear ``runs and tumbles" equation}

The main purpose of the present work is to  provide an alternative approach to study the linear  ``runs and tumbles" (linear RT) equation \eqref{eq:RTeq1}-\eqref{eq:RTeq4} which makes possible to generalize the analysis of \cite{CRS} to any dimension $d \ge 1$. In order to state our main result, we introduce some notations and the functional framework we will work with. 

\smallskip
First, we denote by $m$ some weight function which is either a polynomial or an exponential 
\beqn\label{eq:mPolyExp}
m (x) = \langle x \rangle^k, \ k > 0, 
\quad\hbox{or}\quad
m (x) = \exp ( \gamma \, \langle x \rangle),  \ \gamma \in (0,\gamma^*), 
\eeqn
for some positive constant $\gamma^*$ which will be defined latter, and where $\langle x \rangle^2 = 1 + |x|^2$. 
To a given weight $m$ we define the associated rate function $\Theta_m$ and weight function $\omega$ by 
\bean
\Theta_m(t) 
:= \langle t \rangle^{-\ell}, \ \forall \, \ell \in (0,k), \quad  \omega = 1 \quad &\hbox{if}& \quad m = \langle x \rangle^k;
\\
\Theta_m(t) 
=e^{at}, \Blue \ \forall \, a \in (a^*,0), \Black\quad \omega = m  \quad &\hbox{if}& \quad m = e^{\gamma \, \langle x \rangle},
\eean
\Blue for an optimal rate $a^* = a^*(\gamma) < 0$ which will be also defined later. \Black 
Finally for given weight function $m=m(x) : \R^d \to \R_+$ and exponent $1 \le p \le \infty$, we define the associated weighted Lebesgue space  $L^p(m)$  and weighted Sobolev space $W^{1,p}(m)$,  through their norms 
\be\label{eq:defLpm}
\| f \|_{L^p(m)} := \| m f  \|_{L^p}, \quad  
\| f \|_{W^{1,p}(m)} := \| m f  \|_{W^{1,p}}.
\ee
We use the shorthands $L^p_k = L^p(m)$, when $m = \langle x \rangle^k$, and $H^1(m) = W^{1,2}(m)$. We write $a \lesssim b$ if there exists a positive constant $C$ such that $a \le C \, b$.

\begin{theo}\label{theo:MainTheorem} There exists $\gamma^* > 0$ and there exists a unique positive, invariant by rotations and normalized stationary state 
\bean
&&0 < G \Blue \in L^\infty (m_0) \Black, \quad \langle\!\langle G \rangle\!\rangle = 1, 
\\
&&- v \cdot \nabla_x G + \int_\VV \bigl\{ K' G' - K G \bigr\}  \, dv' = 0,
\eean
where $m_0$ stands for the exponential weight function $m_0(x) := \exp(\gamma_* \langle x \rangle)$. 
Moreover, for any weight function $m$ satisfying \eqref{eq:mPolyExp} and for any $0 \le f_0 \in L^1(m) $, there exists 
a unique solution $f \in C([0,\infty); L^1(m))$ to the  equation  \eqref{eq:RTeq1}-\eqref{eq:RTeq4} associated to
the initial datum $f_0$ and 
\beqn\label{eq:theo2}
\| f(t) - \langle \! \langle f_0 \rangle\! \rangle G \|_{L^1(\omega)} \le \Theta_m(t) \, \|   f_0  \|_{L^1(m)},
\quad \forall \, t \ge 0,
\eeqn
where $\omega$ and $\Theta_m$ are defined just above. 
\end{theo}

\smallskip
The present result generalizes to any dimension $d \ge 1$ similar results (\cite[Theorem~2.1]{CRS} and \cite[Proposition~1]{CRS}) obtained  by Calvez et al.  in dimension $d=1$. 
As pointed out in \cite{CRS}, the main novelty and mathematical interest of the model lie in the fact that the confinement is achieved by a biased velocity jump process, where the bias replaces the confining acceleration field which is the classical confinement mechanism for Boltzmann and Fokker-Planck models, see for instance \cite{MR2562709,MM2016,DMS} and the references quoted therein.  

\smallskip 
\Blue Our strategy is drastically different from the one used in  \cite{CRS} but similar to the approach of \cite{MiSch2016} which develops a Krein-Rutman theory  for positive semigroups which do not fulfill the classical compactness assumption on the associated resolvent but have a nice splitting structure. 
 
However, instead of applying directly the Krein-Rutman Theorem~\cite[Theorem~5.3]{MiSch2016} and in order to be a bit more self-consistent and pedagogical, we rather follow the same line of proof as in \cite{MiSch2016} but we perform some simplifications by  taking advantage of the mass conservation law (or equivalently, that the dual operator has an explicit positive eigenvector). 
The main difficulty is then to get suitable estimates on some related operators and semigroups. 

\smallskip
In section~\ref{sec:L1estim}, the  first step consists in proving a weighted $L^1$ bound which brings out the confinement mechanism. That is the main new bound which is in the spirit of  weighted $L^p$ estimates obtained for performing similar spectral analysis in \cite{GMM*,MiSch2016,MM2016} for Boltzmann, growth-fragmentation and kinetic Fokker-Planck models.  

\smallskip
Next, in order to go further in the analysis, we introduce suitable decompositions 
$$
\LL = A + B, 
$$
with several choices of operators $A$ and $B$ such that $B$ is adequately dissipative, $A$ is $B$-bounded and $A S_B$ enjoys some regularization properties, where $S_B$ stands here for the semigroup associated to the generator $B$. 
For establishing these regularization properties on $AS_B$, we successively use a dispersion argument as introduced by Bardos and Degond \cite{BardosD} for providing better integrability in the position variable (transfer of integrability from the velocity variable to the position variable), next a multiplicator method in the spirit of Lions-Perthame multiplicator (see  \cite{Perthame1,LionsP,Perthame2}) for improving again the integrability estimate in the position variable near the origin and finally a space variable averaging lemma in the spirit of the variant \cite{DesvM,BouchutD} of the classical time and space averaging lemma of Golse et al \cite{GolsePS,GolseLPS}. 
It is worth emphasizing that  the needed regularity estimate is not obtained using an abstract hypocoercivity operator as in \cite{DMS,CRS} nor using an iterated averaging lemma as in  \cite{GMM*} (which allows to transfer regularity from the velocity variable to the position variable thanks to a suitable commutator and the associated ``gliding norms'') but using the more classical (and more robust) averaging lemma. 

\smallskip
More precisely, in Section~\ref{sec:StationaryState}, we make a  first rather simple choice for the splitting of the operator $\LL$ and we obtain that the associated semigroup $S_\LL$ is bounded in the weighted Lebesgue space $X := L^1(m) \cap L^p(m)$ by gathering the above dispersion argument and a shrinkage of the functional space argument as in \cite{MM2016}. A  flavor of the argument reads as follows.  
We write the  iterative version of the Duhamel formula
$$
S_\LL = S_\BB + ... + S_B * (A S_B)^{(*n-1)} + (S_B A)^{(*n)} * S_\LL, \quad \forall \, n \in \N^*, 
$$
where $*$ stands for the usual convolution operator on $\R_+$. We then deduce that $S_\LL$ is bounded in $\BBB(X)$, the space of bounded linear mas from $X$ into itself, by using an exponential decay estimate in $\BBB(X)$  for the terms $S_B * (A S_B )^{(*k)}$, an exponential decay estimate  in $\BBB(L^1,X)$ for $(S_B A)^{(*n)}$ and  exploiting that $S_\LL$ is bounded in  $\BBB(L^1)$ as an immediate consequence of the  mass conservation and the positivity property. 
We then deduce the existence of a weighted uniformly bounded  nonnegative and normalized steady state thanks to a standard Brouwer type fixed point argument. Its uniqueness follows by classical weak and strong maximum principles. As a matter of fact, in the same way one deduces that $0$ is a simple eigenvalue and the only nonnegative eigenvalue. 

\smallskip
In section~\ref{sec:SpectralGap}, we introduce  a more sophisticated surgical truncation $A$ of the kernel operator involved in $\LL$ and thus a second splitting.  
In such a way, using the above mentioned multiplicator method and space averaging lemma, we obtain that the new operator  $A$ is more regular, the corresponding operator $B$ is still  appropriately dissipative for  a suitable equivalent norm and finally $A S_B$ has nice compactness and regularity property. 
With the help of \cite{MiSch2016}, we conclude to a spectral gap on the spectrum of the operator $\LL$ (Weyl theorem) and its translation into an estimate  on the semigroup $S_\LL$ (quantitative spectral mapping theorem) as stated in Theorem~\ref{theo:MainTheorem}.

 \Black

\smallskip
Finally, it is worth pointing out that it is not clear how to use the above analysis in order to make any progress in the understanding of the nonlinear 
equation \eqref{eq:RTeq1}-\eqref{eq:RTeq2}-\eqref{eq:RTeq3}. In particular, we have not been able to prove that the chemical agent density $S$ is decaying with respect to the radial variable $|x|$ and thus that $G$  is also a  stationary state of the nonlinear  equation \eqref{eq:RTeq1}-\eqref{eq:RTeq2}-\eqref{eq:RTeq3}, 
as one can expect by making an analogy with the \Blue one dimension case  and when the velocity set $\VV$  is replaced by  $\bar\VV := \{ -1, 1\}$. \Black Indeed, in that case, 
we may observe  that for $\bar G(x,v) = g(|x|) = C \, e^{-\chi \, |x|}$, we have 
\bean
v \cdot \nabla_x \bar G = v \cdot {x \over |x|} \, g'(|x|) = - \chi \, \hbox{sign} (x\cdot v) \, g(|x|)
=
 \int_{\bar\VV} \bigl\{ K' \bar G - K \bar G \bigr\} \, dv'  , 
\eean 
so that  we have exhibited an explicit (unique, positive and unit mass) stationary state $\bar G$. The associated macroscopic density $\bar \varrho$ is then decaying and thus also the associated  chemical agent density $\bar S$ (thanks to the maximum principle applied to the elliptic equation \eqref{eq:RTeq3}). It turns out then that $\bar G$ is also a  stationary state of the nonlinear  equation \eqref{eq:RTeq1}-\eqref{eq:RTeq2}-\eqref{eq:RTeq3}. 
\Blue We refer the interest reader to the recent paper by Calvez \cite{Calvez*} who establishes the existence of traveling wave solutions for a similar nonlinear model.   \Black

 \medskip
 Let us end the introduction by describing again the plan of the paper. 
 In Section~\ref{sec:L1estim}, we mainly present the weighted $L^1$ estimate which highlights the confinement 
mechanism of the model.    In Section~\ref{sec:StationaryState}, we introduce a first splitting of the generator in order to prove the existence (and next the uniqueness) of a positive stationary state. Finally, in  Section~\ref{sec:SpectralGap}, we introduce a second splitting of the generator which enjoys better smoothness properties
and for which we can use the   Krein-Rutman theory revisited in \cite{MiSch2016} and conclude to the asymptotic stability of the stationary state.

 \medskip
\paragraph{\bf Acknowledgments.} We thank V. Calvez for fruitful discussions which have been a source of the present work. The research leading to this paper was (partially) funded by the French "ANR blanche" project Kibord: ANR-13-BS01-0004.  

\bigskip
\section{Well-posedness and exponential weighted $L^1$ estimate}
\label{sec:L1estim}
\setcounter{equation}{0}
\setcounter{theo}{0}

We first state a well-posedness result concerning the linear RT equation \eqref{eq:RTeq1}-\eqref{eq:RTeq4}, whose very classical proof is skipped,
and we recall some notations and definition. 

\begin{lem}\label{lem:WPSL}
For any $f_0 \in L^p(m)$, $1 \le p \le \infty$, there exists a unique weak (distributional) solution $f \in C([0,T);L^1_{loc}) \cap L^\infty(0,T;L^p(m))$, $\forall \, T \in (0,\infty)$, 
to the linear RT equation \eqref{eq:RTeq1}-\eqref{eq:RTeq4} which furthermore satisfies:

(1) $f(t,\cdot) \ge 0$ for any $t \ge 0$ if $f_0 \ge 0$ (preservation of positivity); 

(2) $\langle\!\langle f(t,\cdot ) \rangle\!\rangle = \langle\!\langle f_0 \rangle\!\rangle $ for any $t \ge 0$ if $L^p(m) \subset L^1$ (mass conservation);  

(3) $f \in C(\R_+;L^p(m))$ and we may associate  to $\LL$ a continuous semigroup $S_\LL$ in $L^p(m)$ by setting $S_\LL(t) f_0 := f(t,\cdot)$ 
for any $t \ge 0$ and $f_0 \in L^p(m)$. Here the continuity has to be understood in the strong (norm) topology sense when $p \in [1,\infty)$ and  in the weak $\sigma(L^\infty,L^1)*$ topology sense when $p=\infty$. 
 \end{lem}

In the sequel, we will thus associate to the generator $\LL$ (and next to the related generators $\BB_0$, $\BB_1$, $\BB$, ...) a semigroup $S_\LL$ and to any initial
datum $f_0 \in L^p(m)$ we will denote by $f(t) = f_\LL(t) = S_\LL(t) f_0$ the solution associated to the related abstract evolution equation.
\Blue
Thanks to by-now standard results (due in particular to Ball~\cite{Ball} and to DiPerna-Lions~\cite{DPL}) this solution is equivalently a distributional, weak or renormalized solution to the corresponding PDE equation. 
\Black
More precisely, the above function $f$ satisfies the linear RT equation \eqref{eq:RTeq1} in the following renormalized sense 
$$
{d \over dt} \int_{\R^d \times \VV} \beta(f) \, \varphi 
=   \int_{\R^d \times \VV} \beta(f) \, v \cdot \nabla_x \varphi +  \int_{\R^d \times \VV}  f \, K \int_\VV [ \varphi' \, \beta(f') - \varphi \, \beta'(f)]  
$$
for any (renormalizing) Lipschitz function $\beta : \R \to \R$ and for any (test) function $\varphi \in C^1_c(\R^d \times \VV)$. It is worth mentioning that 
we deduce the preservation of positivity property by just choosing $\beta(s) = s_-$, $\varphi = 1$, in the above identity and next using the Gronwall lemma.
When $L^p(m) \subset L^1$, the uniqueness result follows from choosing  $\beta(s) = |s|$ and $\varphi = 1$ in the above identity. 
The existence part can be achieved by combining  the characteristics method for the free transport equation and a perturbation by bounded operators argument. 
In the sequel, we will denote $\BBB(X,Y)$ the space of bounded linear operators from a Banach space $X$ into a Banach space $Y$, and we write $\BBB(X) = \BBB(X,X)$. Another way to prove existence and uniqueness consists in using the Hille-Yosida theorem for maximal dissipative unbounded operators. We recall for futur references that we say that an unbounded operator $\LL$ with dense domain $D(\LL) \subset X$ is $a$-dissipative, $a \in \R$, if 
$$
\forall \, f \in D(\LL), \,\, \exists \, f^* \in J_f \quad \langle f^*, \LL f \rangle_{X',X} \le a \, \| f \|_X^2,
$$
where $J_f $ denotes the (nonempty) dual set 
$$
J_f := \{ g \in X'; \,\, \langle g,f \rangle_{X',X} = \| g \|_{X'}^2 = \| f \|_X^2 \}.
$$

\smallskip
We now establish a uniform in time exponential weighted $L^1$ estimate, which is one of the cornerstone arguments of the proof of our main theorem.  

\begin{lem}\label{lem:estimL1m} There exists a constant $\gamma^*> 0$ and for any $\gamma \in (0,\gamma^*)$ there exist  a weight function $ \widetilde{m} $ such that $ \widetilde{m} (x) \sim e^{\gamma \langle x \rangle}$  as $x\to\infty$ and a constant $C \in (0,\infty)$ such that the solution $f$ to the linear RT equation with initial datum $f_0 \in L^1(m)$ satisfies  
\beqn\label{eq:EstimWeightMoment}
 \int |f(t)| \,  \widetilde{m}  \le \max \Bigl( C \!\int |f_0|, \int |f_0| \,  \widetilde{m}  \Bigr), \quad \forall \, t \ge 0.
\eeqn
In particular, the semigroup $S_\LL$ is bounded in $L^1(m)$. 
\end{lem}

\noindent
{\sl Proof of Lemma~\ref{lem:estimL1m}. } We define the dual operator $\LL^*$ by 
$$
\int (\LL^* \varphi) \, f = \int (\LL f) \, \varphi, \quad \forall \, \varphi \in W^{1,\infty}(\R^d \times \VV), \, \, f \in C_c(\R^d \times \VV),
$$
so that 
$$
(\LL^* \varphi)(x,v) := v \cdot \nabla_x \varphi + K  \int_\VV \bigl\{  \varphi' - \varphi \bigr\}  \, dv'.
$$
For a given $\gamma > 0$, we compute 
\bean 
\LL^* \,  e^{\gamma \langle x \rangle} =  \gamma   (v \cdot {x\over \langle x \rangle}) e^{\gamma \langle x \rangle},
\eean
and next
\bean 
\LL^* \Bigl[  (v \cdot {x\over \langle x \rangle}) e^{\gamma \langle x \rangle} \Bigr]
&=&   v \cdot \nabla_x [ (v \cdot {x\over \langle x \rangle}) e^{\gamma \langle x \rangle} ] 
- K \,   (v \cdot {x\over \langle x \rangle}) e^{\gamma \langle x \rangle}
\\
&=&    \Bigl( {|v|^2 \over  \langle x \rangle}  - {(v \cdot x)^2 \over  \langle x \rangle^3}  + \gamma  {(v \cdot x)^2 \over  \langle x \rangle^2} 
 -   (v \cdot {x\over \langle x \rangle})  - \chi \, {|v \cdot x| \over \langle x \rangle} \Bigr) e^{\gamma \langle x \rangle}.
\eean
Defining 
\beqn\label{eq:defV1}
V_1 := \int_\VV|v'_1| \, dv',
\eeqn
\Blue and recalling that  $\zeta$ is  defined in \eqref{eq:RTeq4}, \Black
we finally have
\bean 
\LL^* \Bigl[    {|v \cdot x| \over \langle x \rangle}  \,  e^{\gamma \langle x \rangle} \Bigr]
&=&  
 v \cdot \nabla_x [   {|v \cdot x| \over \langle x \rangle}  \,  e^{\gamma \langle x \rangle} ] 
+  ( 1 + \chi \zeta) (V_1 { |x|   \over \langle x \rangle} -  {|v \cdot x| \over \langle x \rangle} )  e^{\gamma \langle x \rangle}
\\
&=&  \zeta \, \Bigl(  { |v|^2  \over \langle x \rangle} - {|v\cdot x|^2 \over \langle x \rangle^3} 
+ \gamma \,{|v\cdot x|^2 \over \langle x \rangle^2}  \Bigr)  \,  e^{\gamma \langle x \rangle}
\\
&&+  \Bigl(   \,  ( 1 + \chi \zeta) V_1 { |x|   \over \langle x \rangle} -  ( 1 + \chi \zeta)  {|v \cdot x| \over \langle x \rangle} \Bigr)  \,  e^{\gamma \langle x \rangle}.
 \eean
\Blue Defining the weight function \Black
\beqn\label{eq:def:weight:mL1}
 \widetilde{m} :=  \Bigl(  1 + \gamma (v \cdot {x \over \langle x \rangle}) - \beta \,  {|v \cdot x| \over \langle x \rangle}  \Bigr) e^{\gamma \langle x \rangle}, 
\eeqn
\Blue for  $\beta,\gamma \in (0,1)$  to be fixed precisely later,  we observe that  $\tilde m$ satisfies 
\beqn\label{eq:def:weight:mL1bis}
 \exists \, \delta \in (0,1), \quad e^{\gamma \langle x \rangle} (1-\delta)  \le  \widetilde{m}  \le (1+\delta) \,e^{\gamma \langle x \rangle},  
\eeqn
 by choosing $\beta,\gamma $ small enough and because $\VV := B(0,V_0)$ is a bounded set. \Black Gathering the previous identities, we find
\bean
(\LL^*  \widetilde{m}  ) \, e^{-\gamma \langle x \rangle}
&=& 
\gamma  \Bigl( {|v|^2 \over  \langle x \rangle}  - {(v \cdot x)^2 \over  \langle x \rangle^3}  + \gamma  {(v \cdot x)^2 \over  \langle x \rangle^2} 
 - \chi \, {|v \cdot x| \over \langle x \rangle} \Bigr)  
\\
&&
 - \beta   \, \zeta \, \Bigl(  { |v|^2  \over \langle x \rangle} - {|v\cdot x|^2 \over \langle x \rangle^3} 
+ \gamma \,{|v\cdot x|^2 \over \langle x \rangle^2}  \Bigr)   
 - \beta  
  \Bigl(   \,  ( 1 + \chi \zeta) V_1 { |x|   \over \langle x \rangle} -  ( 1 + \chi \zeta)  {|v \cdot x| \over \langle x \rangle} \Bigr) 
\\
&\le& 
\gamma    V_0^2  \Bigl( {1 \over \langle x \rangle}  + \gamma \Bigr)  
 +\beta   V_0^2  \Bigl( {2 \over \langle x \rangle} + \gamma   \Bigr)  
+\beta (1+\chi)   {|v \cdot x| \over \langle x \rangle}  
\\
 &&  
- \gamma \chi \,     {|v \cdot x| \over \langle x \rangle}  
- \beta (1-\chi) V_1   {|x| \over \langle x \rangle}  
\\
&\le& 
[V_0^2  (\gamma + 2 \beta) + \beta (1-\chi) V_1] {1  \over \langle x \rangle}    
+ [\gamma^2V_0^2 + \beta\gamma V_0^2 -  \beta (1-\chi) V_1 ]. 
\eean
We thus deduce that 
\beqn\label{eq:L*weight:mL1}
(\LL^*  \widetilde{m}  ) \, e^{-\gamma \langle x \rangle}
\le
 {C \over \langle x \rangle} - 2 \alpha ,
\eeqn
by choosing $\beta (1+\chi) = \gamma \chi$ and $\gamma > 0$ small enough in such a way  that $2\alpha :=  \beta (1-\chi) V_1 - \gamma^2V_0^2  - \beta\gamma V_0^2 > 0$.  
Observing that 
$$
\Bigl( {C \over \langle x \rangle} - 2 \alpha \Big) \, e^{\gamma \langle x \rangle} 
\le C e^{\gamma \langle R \rangle} \, {\bf 1}_{B(0,R)} - \alpha (1+\delta) \, e^{\gamma \langle x \rangle} 
\le A  - \alpha   \widetilde{m} ,
$$
for some constant $A > 0$, we have proved
$$
\LL^*  \widetilde{m}  \le A - \alpha \,  \widetilde{m} , \quad \alpha < 0.
$$

\smallskip
We consider now $f$ the solution to the linear RT equation \eqref{eq:RTeq1} associated to $f_0 \in L^1(m)$. Denoting  $g := |f|$, we deduce from the above inequality
that 
\bean
{d \over dt} \int g \,  \widetilde{m}  
=  \int (\LL f) \, (\hbox{sign} f ) \,  \widetilde{m}  
\le  \int (\LL g)  \,  \widetilde{m} 
\le A \int g   - \alpha  \int g   \widetilde{m} . 
 \eean
As a consequence, \eqref{eq:EstimWeightMoment} holds with $C := A/\alpha$. \qed

\bigskip
\section{The stationary state problem}
\label{sec:StationaryState}
\setcounter{equation}{0}
\setcounter{theo}{0}

We introduce two generators $\BB_0$ and next $\BB_1$  in the following paragraphs and we study the
gain of integrabilty properties of the associated semigroups. \Blue Introducing the splitting $\LL = \AA_1 + \BB_1$, we then use these estimates in order to prove that $S_\LL$ is a bounded semigroup in $L^1(m) \cap L^p(m)$, from what we deduce the existence of a stationary state. Uniqueness of this one is finally proved thanks to classical weak and strong maximum principles. 
\Black

\subsection{The operator $\BB_0$ and the associated semigroup $S_{\BB_0}$. } 
We define $\BB_0$ by 
$$
\BB_0 f := - v \cdot \nabla_x f - K f . 
$$

\begin{lem}\label{lem:SB0Lpm} 
There exist $\gamma^* > 0$ and $a^* < 0$ such that for any  $1 \le p \le \infty$,   $m = e^{\gamma \langle x \rangle}$, $\gamma \in [0,\gamma^*)$, there holds  
$$
\| S_{\BB_0} (t) \|_{L^p(m) \to L^p(m)} \lesssim e^{ a t }, \quad \forall \, a > a^*.
$$
\end{lem}

\noindent
{\sl Proof of Lemma~\ref{lem:SB0Lpm}. } We consider a solution $f = S_{\BB_0} (t) f_0$ to the evolution equation 
associated to $\BB_0$ and we compute
\bean
{d \over dt} \int |f|^p m^p 
&=& \int  |f|^p \,p \,  [ v \cdot \nabla_x m - K \, m]  \, m^{p-1}
\\
&\le& p \, [V_0 \, \gamma + \chi - 1] \int   |f|^p m^p ,
\eean
from which we easily conclude thanks to the Gronwall lemma. 
\qed

\medskip
For any $\varphi \in L^\infty_{xv}$, we define $A = A_\varphi : L^q_xL^1_v \to L^q_x$ by 
$$
A f = A_\varphi f (x): = \int_\VV \varphi(x,v') \, f(x,v') \, dv'.
$$
Given some Banach spaces $X_i$ and two functions $u \in L^1(\R_+;\BBB(X_2,X_3))$, $v \in L^1(\R_+;\BBB(X_1,X_2))$, we define the convolution function
$u * v \in L^1(\R_+;\BBB(X_1,X_3))$ by 
$$
(u*v)(t) := \int_0^t u(t-s) \, v(s) \, ds.
$$
We also define $u^{*n}$ by $u^{*1} = u$ and $u^{*n} = u^{*(n-1)}*u$ for $n \ge 2$.

\begin{lem}\label{lem:SB0L1LptoLp} \Blue There exists $\gamma^* > 0$ and for any $\gamma \in [0,\gamma^*)$ there exists $a^* < 0$ such that for $m = e^{\gamma \langle x \rangle}$ and for any $\varphi \in L^\infty_{xv}$,  
there holds
\beqn\label{eq:Reg2S0}
\| A_\varphi S_{\BB_0} (t) \|_{L^1_xL^\infty_v(m) \to L^\infty_{xv}(m)  } \lesssim  t^{-d} \,e^{a t}, \quad \forall \, t > 0, \,\,  \forall \, a > a^*.
\eeqn
As a consequence,  there exists $n \in \N^*$ ($n = d+2$ is suitable) such that  
\beqn\label{eq:Reg2S0TER}
\| (A_\varphi S_{\BB_0})^{(*n)} (t) \|_{L^1_{xv}(m) \to L^\infty_{xv}(m) } \lesssim   e^{a t}, \quad \forall \, t \ge  0, \,\,  \forall \, a > a^*.
\eeqn
 
\end{lem}

\noindent
{\sl Proof of Lemma~\ref{lem:SB0L1LptoLp}. } We split the proof into two steps.

\smallskip\noindent
{\sl Step 1. } We adapt the classical  dispersion result of Bardos and Degond \cite{BardosD}. 
We denote $f(t) := S_{\BB_0}(t) f_0$ the solution to the damped transport equation 
$$
\partial_t f+ v \cdot \nabla_x  f+ K f  = 0, \quad f(0) = f_0.
$$
The characteristics method gives the representation formula
$$
(S_{\BB_0}(t) f_0 ) (x,v) = f_0(x-vt,v) \, e^{-\int_0^t K(x-vs,v) ds}. 
$$
We then have 
$$
\rho(t,x) := (A S_{\BB_0} (t) f_0) (x) = \int_\VV \varphi(x,v_*) f_0(x-v_*t,v_*) \, e^{-\int_0^t K(x-v_*s,v_*) ds} \, dv_*. 
$$
\Blue
Using that $K(x,v) \ge 1-\chi$ and $\langle x \rangle \le \langle x - v_* t\rangle + |v_*| t + 1$, we deduce 
\bean
| \rho(t,x) |  
&\le&   e^{-(1-\chi) t} \int_\VV | \varphi(x,v_*) |  \bigl[ \sup_{w \in \VV} |f_0|(x-v_* t, w) \bigr] e^{\gamma \langle x - v_* t\rangle} \, dv_* \, e^{\gamma (V_0 t +1 - \langle x \rangle) }
\\
&\le& e^{(\chi + \gamma V_0 - 1)  t + \gamma} \| \varphi  \|_{L^\infty}  \int_{\R^d}   \bigl[ \sup_{w \in \VV} |f_0|(y, w) \bigr] e^{\gamma \langle y \rangle}  \, {d y \over t^d} \, e^{- \gamma \langle x \rangle}.
\eean
Defining $\gamma^* := (1-\chi)/V_0$ and for any $\gamma \in [0,\gamma^*)$ defining $a^* := \chi + \gamma V_0 - 1 < 0$, we conclude  with
$$
\|  A S_{\BB_0}(t) f_0 \|_{L^\infty_{xv}(m)} \lesssim {e^{a^* t}  \over t^d} \, \| f_0 \|_{L^1_xL^\infty_v(m)},
$$
which in particular implies \eqref{eq:Reg2S0}.  \Black

%

\smallskip\noindent
{\sl Step 2. } From Lemma~\ref{lem:SB0Lpm},  for $r = 1$ and $r=\infty$, we clearly have
\bean
\|  A S_{\BB_0}(t)  \|_{\Blue L^r_{x}L^\infty_v(m) \to L^r_{x}L^\infty_v(m)} \lesssim
e^{ a \, t }.
\eean
Gathering that estimate with \eqref{eq:Reg2S0}, we may repeat the proof of \cite[Proposition 2.2]{MQT2016} (see also \cite[Lemma 2.4]{MM2016} or \cite{Mbook*}) and we get 
\beqn\label{eq:Reg2S0eq3}
\|  (A S_{\BB_0})^{(*d+1)} (t)  \|_{\Blue L^1_{x}L^\infty_v(m) \to L^\infty_{x}L^\infty_v(m)} \lesssim
e^{ a \, t }.
\eeqn
 Observing that  
\beqn\label{eq:Reg2S0eq2}
\|  (A S_{\BB_0})  (t)  \|_{L^1_{x}L^1_v(m) \to L^1_{x}L^\infty_v(m) } \lesssim
e^{ a \, t },
\eeqn
thanks to Lemma~\ref{lem:SB0Lpm} and $A : L^1_{xv} (m) \to L^1_xL^\infty_v(m)$, we  conclude to \eqref{eq:Reg2S0TER} by taking $n = d+2$. \qed

\subsection{The operator $\BB_1$ and the associated semigroup $S_{\BB_1}$. } \label{subsect:B1}
We define $\BB_1$ by 
\beqn\label{eq:defBB1}
\BB_1 f := - v \cdot \nabla_x f - K f + (1-\phi_R) \int_\VV K' f' \, dv' ,
\eeqn
 where we have defined the truncation functions  $\phi_R (x) := \phi(x/R)$ for a given 
  $\phi \in \DD(\R^d)$ which is radially symmetric and satisfies ${\bf 1}_{B(0,1)} \le \phi \le {\bf 1}_{B(0,2)}$. 

\begin{lem}\label{lem:SB1L1m} There exist $\gamma^* > 0$, $a^* < 0$ and $C \ge 1$ such that for $R$ large enough and $m = e^{\gamma \langle x \rangle}$, $\gamma \in (0,\gamma^*)$, the semigroup $S_{\BB_1}$ satisfies the following
growth estimate
$$
\| S_{\BB_1} (t) \|_{L^1(m) \to L^1(m)} \le C \, e^{a^*t}, \quad \forall \, t \ge 0. 
$$
\end{lem}

\noindent
{\sl Proof of Lemma~\ref{lem:SB1L1m}. }
We observe that the dual operator $\BB_1^*$ writes
$$
\BB_1^* \varphi = \LL^* \varphi - \phi_R \, K  \varphi.
$$
Defining the modified weight function $\widetilde{m}$ as in \eqref{eq:def:weight:mL1} and using the inequalities \eqref{eq:def:weight:mL1bis} and  \eqref{eq:L*weight:mL1},
we get 
\bean
\BB_1^* \widetilde{m} 
&\le& 
\Bigl(  {C \over \langle x \rangle} - 2 \alpha \Bigr)  e^{-\gamma \langle x \rangle} - \phi_R \, K \,  \widetilde{m} 
\\
&\le& 
\Bigl(  {C \over \langle x \rangle} - 2 \alpha \Bigr)  e^{-\gamma \langle x \rangle} - {\bf 1}_{B_{R}(0)} \, (1-\chi) \, (1-\delta) \, e^{-\gamma \langle x \rangle}
\\
&\le&  -  \alpha \, \widetilde{m},
\eean
for $\gamma> 0$ small enough and $R> 1$ large enough, because $C = \OO(\gamma)$ and $\delta = \OO(\gamma)$. We then have proved that 
$\BB_1$ is dissipative in $L^1(\widetilde{m})$ and we immediately conclude. \qed
 
\begin{lem}\label{lem:SB1LpLp} For the same constants $\gamma^* > 0$ and $a^* < 0$ as defined in Lemma~\ref{lem:SB1L1m},   for any $1 \le p \le \infty$ and $m = e^{\gamma \langle x \rangle}$, $\gamma \in (0,\gamma^*)$,  the semigroup $S_{\BB_1}$ satisfies the growth estimate 
\beqn\label{eq:SB1LpLp}
\| S_{\BB_1} (t) \|_{\Blue L^1(m) \cap L^p(m) \to L^1(m) \cap L^p(m)} \le C \, e^{a t}, \quad \forall \, t \ge 0, \quad \forall \, a > a^*.
\eeqn
\end{lem}

\noindent
{\sl Proof of Lemma~\ref{lem:SB1LpLp}. } We write 
$$
\BB_1 = \BB_0 + \AA^c_0,
$$
with  $ \AA^c_0 = A_\psi$, $\psi:= (1-\phi_R) K(x,v)$, and then the iterated Duhamel formula 
$$
S_{\BB_1} = \bigl\{ S_{\BB_0} + ... + S_{\BB_0} * (\AA^c_0  S_{\BB_0})^{(*n)}  \bigr\}
+ S_{\BB_0} * (\AA^c_0  S_{\BB_0})^{(*n)}*\AA^c_0S_{\BB_1}   =: U_1 + U_2,
$$
with $n = d+2$.  From Lemma~\ref{lem:SB0Lpm}, \eqref{eq:Reg2S0TER}
in Lemma~\ref{lem:SB0L1LptoLp} and Lemma~\ref{lem:SB1L1m}, we easily deduce
$$
\| U_2 \|_{L^1\to L^\infty} \le \| S_{\BB_0} \|_{L^\infty \to L^\infty} * \| (\AA^c_0  S_{\BB_0})^{(*n)}\|_{L^1 \to L^\infty}*\|\AA^c_0S_{\BB_1}\|_{L^1\to L^1}   
\lesssim e^{a t}, 
$$
 for any $a > a^*$, where we have removed the weight dependence to shorten notations. We have similarly the same decay estimate on the remainder term $U_1$ in $\BBB(X)$ by just using Lemma~\ref{lem:SB0Lpm}. We
deduce that \eqref{eq:SB1LpLp} holds for $p=\infty$. We conclude that  \eqref{eq:SB1LpLp} holds for any $1 \le p \le \infty$ by interpolating that first estimate in $L^\infty$  together with the
estimate established in Lemma~\ref{lem:SB1L1m}. \qed

\begin{lem}\label{lem:SB1L1LptoLp} 
For the same constants $\gamma^* > 0$ and $a^* < 0$ as defined in Lemma~\ref{lem:SB1L1m}, 
 for any $\varphi \in L^\infty_{xv}$ and $m = e^{\gamma \langle x \rangle}$, $\gamma \in [0,\gamma^*)$,  there holds
\beqn\label{eq:SB1L1LptoLp}
\| A_\varphi S_{\BB_1} (t) \|_{\Blue L^1_x L^\infty_v(m) \to L^\infty_{xv}(m)}  \lesssim  t^{-d}   e^{a t}, \quad \forall \, t \ge 0, \quad \forall \, a > a^*.
\eeqn
As a consequence, for 
 $n \in \N^*$  large enough ($n= d+2$ is suitable), there holds
\beqn\label{eq:ASB1L1toLinfty}
\| (A_\varphi S_{\BB_1})^{(*n)} (t) \|_{\Blue L^1(m) \to  L^\infty(m)} \lesssim  e^{a t}, \quad \forall \, t \ge 0, \quad \forall \, a > a^*.
\eeqn
\end{lem}

\noindent
{\sl Proof of Lemma~\ref{lem:SB1L1LptoLp}. } With the notation of the proof of Lemma~\ref{lem:SB1LpLp}, we write  
$$
A_\varphi S_{\BB_1} =   A_\varphi S_{\BB_0} + A_\varphi S_{\BB_0}*  \AA^c_0  S_{\BB_1},
$$
and we immediately conclude that \eqref{eq:SB1L1LptoLp} holds putting together \eqref{eq:Reg2S0} and the fact that $\AA^c_0  S_{\BB_1}$ has the 
appropriate decay rate in $\BBB(L^1(m);L^1_x L^\infty_v(m))$ thanks to Lemma~\ref{lem:SB1L1m}. Introducing the notations $X_1 := L^1_x L^\infty_v(m)$
and $X_\infty := L^1_x L^\infty_v(m) \cap L^\infty_{xv}(m)$, we then easily see from \eqref{eq:SB1LpLp} and \eqref{eq:SB1L1LptoLp} that 
$$
\| A_\varphi S_{\BB_1} \|_{X_p \to X_q} \lesssim \Theta_{p,q}(t) \, e^{at}, \quad \forall \, a > a*,
$$
for $(p,q) = (1,1)$, $(1,\infty)$, $(\infty,\infty)$, with $\Theta_{1,1} = \Theta_{\infty,\infty} = 1$ and $\Theta_{1,\infty} = t^{-d}$.  
Repeating again the proof of \cite[Proposition 2.2]{MQT2016},  we deduce 
\bean
\|  (A_\varphi S_{\BB_1})^{(*n-1)} (t)  \|_{X_1 \to X_\infty } \lesssim
e^{ a \, t }.
\eean
We conclude by using  that $A_\varphi  S_{\BB_1}$ has the 
appropriate decay rate in $\BBB(L^1(m);L^1_x L^\infty_v(m))$ thanks to Lemma~\ref{lem:SB1L1m}.  \qed

\subsection{Existence of a steady state} We establish now the existence of a steady state thanks to a fixed point argument. 
We fix an exponential weight $m := e^{\gamma \langle x \rangle}$, $\gamma \in (0,\gamma^*)$. 
We define the Banach space {\Blue$X := L^1(m) \cap L^\infty(m)$} as well as 
$$
\forall \, f \in X, \quad \Nt f \Nt := \sup_{t \ge 0}  \|  S_\LL(t) f \|_X. 
$$

\begin{lem}\label{lem:FixedPt}  \Blue The semigroup $S_\LL$ is bounded in $X$. As a consequence, there exists at least one nonnegative, invariant by rotation and normalized stationary state $G \in X$ to the linear RT equation \eqref{eq:RTeq1}. \end{lem}

\noindent
{\sl Proof of Lemma~\ref{lem:FixedPt}. }  We split the proof into two steps. 

\smallskip\noindent{\sl Step 1.} 
We split the operator $\LL$ as 
$$
\LL = \AA_1 + \BB_1, \quad \AA_1  : = A_\psi, \quad \psi:= \phi_R K(x,v),
$$
and with the same integer $n$ as in Lemma~\ref{lem:SB1L1LptoLp} we write the iterated Duhamel formula
$$
 S_\LL = \bigl\{ S_{\BB_1} + ... + S_{\BB_1} * (\AA_1 S_{\BB_1})^{(n*)}  \bigr\} +
 S_{\BB_1} * (\AA_1 S_{\BB_1})^{(n*)}  * \AA_1 S_\LL =: V_1 + V_2.
$$
For the first term and thanks to Lemma~\ref{lem:SB1LpLp}, for some constant $K_1 \ge 1$,  we have 
$$
\| V_1(t)  f_0 \|_{X}  \le \sum_{\ell=0}^n\| S_{\BB_1} *  (\AA_1 S_{\BB_1})^{(\ell *)} f_0  \|_{X}  \le  K_1 \,   \| f_0 \|_{X} .
 $$
For the second term, for some constant $K_2 \ge 1$, we have 
$$
\| V_2(t)   \|_{L^1(m) \to X} \le \|  S_{\BB_1} \|_{X \to X} * \| (\AA_1 S_{\BB_1})^{(n*)} \|_{L^1(m) \to X}  * \| \AA_1 S_\LL \|_{L^1(m) \to
L^1(m) } \le K_2,
$$
using both \eqref{eq:SB1LpLp} and \eqref{eq:ASB1L1toLinfty}. 
  All together, we have proved
\bean
\| S_\LL(t)  f_0 \|_{X}  \lesssim    \| f_0 \|_{X}, \quad \forall \, t \ge 0.
\eean
As a consequence, the quantity $\Nt \cdot \Nt$ defines a norm on  $X$ which is equivalent to its usual norm $\| \cdot \|_X$.

\smallskip\noindent{\sl Step 2.} 
\Blue
For a given $g_0 \in X$  which is also invariant by rotation and a probability measure, we define $C := \Nt g_0 \Nt$ and next  the set 
$$
\CC := \Bigl\{ f \in   X; \,\, f \ge 0, \,\, \langle\!\langle f \rangle\!\rangle = 1, \,\, f_R = f, \, \forall \, R \in SO(d), \,\, \Nt f \Nt  \le C \Bigr\}, 
$$
which is not empty (e.g. $g_0 \in \CC$), convex and compact for the weak $*$ topology of $X$. 
Moreover, thanks to Lemma~\ref{lem:WPSL}, the flow is continuous for the $L^1$ norm and preserves positivity and total mass. By construction, we see that for any $f_0 \in \CC$ and  $t \ge 0$, we have 
$$
\Nt S_\LL(t)  f_0 \Nt  = \sup_{s \ge t} \| S_\LL(s)  f_0 \|_{L^p(m)} \le \sup_{s \ge 0} \| S_\LL(s)  f_0 \|_{L^p(m)} = \Nt  f_0 \Nt \le C.
$$
All together, the set $\CC$ is clearly  invariant by the flow $S_\LL$. Thanks to a standard variant of the Brouwer-Schauder-Tychonoff fixed point theorem
(see for instance \cite[Theorem 1.2]{EMRR}), we obtain the existence of an invariant element $G$ for the linear RT flow which furthermore belongs to $\CC$, in other words
$$
\exists \, G \in \CC \,\,\hbox{ such that } \,\,  S_\LL(t) G = G, \,\, \forall \, t \ge 0. 
$$
As a consequence, we have $G \in D(\LL) \backslash \{  0 \}$ and $\LL G = 0$, so that  $G$ is a stationary state for the  linear RT equation which fulfills all the properties listed in the statement of Theorem~\ref{theo:MainTheorem}. 
\qed
\Black


\subsection{Uniqueness of the stationary state}

In this section we prove  a weak and a strong maximal principle on the operator $-\LL$. The uniqueness of the normalized and positive steady state then follows using classical arguments. We skip the proof of that last one and we refer for instance to \cite[Step 4, proof of Theorem~5.3]{MiSch2016} for details. 
 
 \Blue
 \begin{lem}\label{lem:*MP} 
 The operator $\LL$ satisfies the following Kato's inequality
\bear\label{eq:KatoIneq}
(\hbox{sign}f) \LL f := {1 \over 2|f|} ( f \LL \bar f + \bar f \LL f) &\le& \LL|f|,
\eear
for any complex valued function $f \in X + i X$. As a consequence, the operator $-\LL$ satisfies the weak maximum principle. 
\end{lem}

\noindent{\sl Proof of Lemma~\ref{lem:*MP}. } For $f \in X + i X$, $f \not= 0$, we just compute 
\bean 
{1 \over 2|f|} ( f \LL \bar f + \bar f \LL f) 
&=& -v\cdot\nabla|f|-K|f  + {1 \over 2|f|} \Bigl( f \int_\VV K' \bar f' + \bar f \int_\VV K' f' \Bigr)
\\
&\le& -v\cdot\nabla|f|-K|f | +\int_\VV K'|f'|=\LL|f|.
\eean
The weak maximum principle follows from Kato's inequality \eqref{eq:KatoIneq} and the mass conseravtion $\LL^* 1 = 0$  thanks to 
 classical  characterizations of positive semigroups in \cite{Arendt,Schep}.
\qed

\Black

\begin{lem}\label{lem:-Lspos}
The operator $-\LL$ satisfies the following version of the strong maximum
principle: for any given $0 \le g\in L^2(m)$ and $\lambda\in\R$, there holds
$$
g\in D(\LL)\setminus\{0\} \ and \ (-\LL+\lambda)g\geq0 \ imply \ g>0
$$
\end{lem}

\noindent{\sl Proof of Lemma~\ref{lem:-Lspos}. } We consider $g$ as in the above statement and we prove that it is 
a positive function in several steps. 

\smallskip\noindent
{\sl Step 1. } Defining $M := 1 +\chi + \lambda \in \R$ and $m := 1-\chi > 0$, we see that 
\bean
v \cdot \nabla_x g + M \, g 
&\ge& v \cdot \nabla_x g +K \, g \Blue + \lambda \, g 
\\
&\ge& \int_\VV K' \, g' \, dv' 
\ge   m \, \varrho, \quad  \varrho := \int_\VV  g' \, dv'. 
\eean
Integrating the above inequality along the free transport characteristics, we get 
\beqn\label{eq:StrongPMStep1}
g(x,v) \ge  m \int_0^t \varrho (x - v \, s) \, e^{- M \, s} \, ds 
+ g (x- v \, t,v) \, e^{- M \, t}, \quad \forall \, t \ge 0.
\eeqn

\smallskip
\noindent {\sl Step 2. } In particular, because the second term at the RHS is nonnegative, we may keep only the contribution of the first  term  and we get
\bear\nonumber
\varrho(x)  &\ge&  m \int_{1/2}^1 \int_\VV \varrho (x - v \, s) \, e^{- M} \, ds dv
\\  \nonumber 
&\ge& \kappa \int_{\VV/2} \varrho (x + w) dw, \quad \kappa > 0.
\eear
Since $g \ge 0$ and $g \not\equiv 0$, we also have $\varrho \ge 0$ and $\varrho \not\equiv 0$,  and  there exists $x_0 \in \R^d$ and $r > 0$, small enough, such that 
$\langle \varrho \, {\bf 1}_{B(x_0,r)} \rangle = \alpha > 0$ and $B(0,2r) \subset \VV$, or in other words $B(x_0,r) \subset x + \VV/2$ for any $x \in B(x_0,r)$. 
As a consequence, we have 
\beqn\label{eq:StrongPMStep2}
\varrho \ge \alpha_0 \, {\bf 1}_{B(x_0,r)}, \quad \alpha_0 := \kappa \, \alpha.
\eeqn

\smallskip
\noindent {\sl Step 3. } Observing that for any $x \in \R^d$, there exists a small ball $B \subset \VV$ and some times $\tau_1 > \tau_0 > 0$ such that $x - sv \in B(x_0,r/2)$ for any $v \in B$ and $s \in (\tau_0,\tau_1)$, we may argue as above and we get 
\bear\nonumber
\varrho(x)  &\ge&  m \int_{\tau_0}^{\tau_1} \int_B \varrho (x - v \, s) \, e^{- Ms} \, ds dv
\\  \nonumber 
&\ge& \alpha_x := m |B| (\tau_1-\tau_0) \, e^{-M \tau_1} \, \alpha_0  > 0.
\eear
Finally, using \eqref{eq:StrongPMStep1} again, we deduce 
$$
g(x,v) \ge   m \int_0^1 \varrho (x - v \, s) \, e^{- M \, s} \, ds  > 0, 
$$
which  concludes the proof.   \qed

%
%

\bigskip
\section{Asymptotic stability of the stationary state}
\label{sec:SpectralGap}
\setcounter{equation}{0}
\setcounter{theo}{0}

\subsection{A new splitting }
In all this section, excepted in paragraph~\ref{subseq:ASEinL1}, we fix an exponential  weight function $m  = e^{\gamma \langle x \rangle}$, with $\gamma \in (0,\gamma^*)$ and $\gamma^* > 0$ defined in Lemma~\ref{lem:SB1L1m}, and we define the Banach space {\Blue$X := L^1(m) \cap L^2(m)$}. 
We also introduce the second splitting
\beqn\label{eq:defL=A+B}
\LL = \AA + \BB, \quad \AA f :=  \int_\VV K'_{R,\delta_i} f' \, dv',
\eeqn
where 
$$
 \quad K_{R,\delta_i} = \phi_{\delta_2,R}(x) \psi_{\delta_1} (v)\, K_{\delta_3} (x,v), \quad K_{\delta_3} (x,v) = 1 + \chi \zeta_{\delta_3}(x \cdot v), 
$$
for some real numbers $R > 1$, $\delta_1,\delta_2, \delta_3 \in (0,1)$ to be fixed, 
and where we have defined the truncation functions  $\phi_\lambda(z) := \phi(z/\lambda)$,   $\phi \in \DD(\R^d)$ radially symmetric, ${\bf 1}_{B(0,1)} \le \phi \le {\bf 1}_{B(0,2)}$, and then  $\phi_{\delta,R} (x) := \phi_R (x) - \phi_\delta (x)$,  $\psi_{\delta} (v) := 1 - \phi_\delta (V_0 - |v|)  - \phi_\delta (v)$, as well as a regularized sign function $\zeta_\delta \in C^\infty(\R)$ which is odd, increasing and satisfies $\zeta_\delta(s) = 1$ for any $s \ge \delta$. 
\Blue It is worth mentioning that the kernel $K_{R,\delta_i} \in C^\infty_c(\R^d \times \R^d)$ with support included in $B(0,2R) \cap B(0,V'_0)$, $V'_0 \in (V_0-\delta_1,V_0)$, what will be a cornerstone used property during the proofs of  Proposition~\ref{prop:SBLpLp} and  Proposition~\ref{prop:Reg3}. That smoothness property contracts with the non-smoothness property of the kernel $\phi_R K$ associated to the operator $\AA_1$ and that is the main reason for introducing that new splitting.\Black
 
\smallskip
We establish that $S_\BB$ and $\AA S_\BB$ enjoy suitable decay estimate (Section~\ref{subseq:DEstimSB}) and regularity estimate 
(Section~\ref{subseq:REstimASB}) from which we deduce the asymptotic stability in $X$ thanks to a semigroup version of the Krein-Rutman theorem 
(Section~\ref{subsec:ASEinX}) and next the asymptotic stability in any exponential and polynomial weighted $L^1$ space by using an extension argument (Section~\ref{subseq:ASEinL1}).

\subsection{Decay estimates for the semigroup $S_{\BB}$ }
\label{subseq:DEstimSB}

\begin{prop}\label{prop:SBLpLp} For the same constant  $a^* < 0$ as defined in Lemma~\ref{lem:SB1L1m},  there holds
\beqn\label{eq:SBX}
\| S_{\BB} (t) \|_{X \to X} \lesssim e^{a t}, \quad \forall \, t \ge 0, \quad \forall \, a > a^*.
\eeqn
\end{prop}

\noindent
{\sl Proof of Proposition~\ref{prop:SBLpLp}. } We split the proof into five steps. 

\smallskip\noindent
{\sl Step 1. Norms and splitting. } Inspired by \cite[Proposition~5.15]{GMM*} and the moment trick introduced in \cite{LionsP}, we define the three norms $\| \cdot \|_X$, $\Nt \cdot \Nt$ and $N(\cdot)$ in the following way
\bean
\| f \|_X^2 &:=& \| f \|_{L^1(m)}^2 + \| f \|_{L^2(m)}^2, 
\\
\Nt f \Nt^2 &:=& \eta_2 \| f \|^2_X + \int_0^\infty \| S_{\BB_1}(\tau) f \|^2_X \, d\tau, 
\\
N(f)^2 &:=& \eta_1 \, \| f \|_{L^2(\mu^{1/2})}^2 + \Nt f \Nt^2, 
\eean
for some constants $\eta_1, \eta_2 \in (0,1)$ to be fixed and where $\mu$ is the weight function
$$
\Blue\mu :=  \Bigl( 1 - {x \over |x|^{1/2}} \cdot {v \over |v|} \Bigr) \, \phi_{1/2}(x), 
$$
so that $0 \le \mu \le 2 \, \phi_1$. 
Thanks to the decay estimate of Lemma~\ref{lem:SB1LpLp}, one easily sees that these three norms are equivalent. 
\Blue In particular, there exists $a^* := a^*(m,\eta_i) < 0$ such that 
\beqn\label{eq:equivNX}
\forall \, f \in X, \quad - \| f \|_X \le 4 a^* N(f).
\eeqn\Black

\smallskip
 We fix $f_0 \in X$ and we define $f(t) = f_\BB(t) = S_\BB(t) f_0$ the associated trajectory along the action of the semigroup $S_\BB$. 
\Blue Our goal is to establish  \eqref{eq:SBX} by proving that $\BB$ is suitably dissipative for the norm $N(\cdot)$. \Black In order to do so, we write
\beqn\label{eq:defT=T123}
T := {1 \over 2} {d \over dt} N( f_\BB (t) )^2  =  \eta_1 T_1 + \eta_2  T_2 + T_3,
\eeqn
where $T_i$ are the contributions of the terms involved in the definition of the norm $N(\cdot)$ that we compute separately. For latter references, we introduce the splitting 
of  $\BB$  as 
\bean
\BB  
&:=&\BB_1 + \AA^c_1 + \AA^c_2 + \AA^c_3,
\eean
where $\BB_1$ is defined in section~\ref{subsect:B1} with $R \ge 1$ large enough so that Lemma~\ref{lem:SB1L1m} and Lemma~\ref{lem:SB1LpLp} hold true, and we have set 
\bean
\AA^c_1 f &=&  \phi_R(x)  \int_\VV K' f' \psi^c_{\delta_1}(v') \, dv'
\\
\AA^c_2 f &=&  \phi_{\delta_2}(x)  \int_\VV K' f' \psi_{\delta_1}(v')  \, dv'
\\
\AA^c_3 f &=& \phi_{\delta_2,R}(x)  \int_\VV K_{\delta_3}^c(x \cdot v') f' \psi_{\delta_1}(v') \, dv',
\eean
with  $ \psi^c_{\delta_1} := 1 -\psi_{\delta_1}$, $\phi^c_R := 1 - \phi_R$,  $K_\delta^c =  \chi \zeta^c_\delta$, $ \zeta^c_\delta = \zeta - \zeta_\delta$. 
We shall also denote
$$
\AA^c_{123} := \AA^c_1 + \AA^c_1 + \AA^c_2
\quad\hbox{and}\quad
\AA^c_{0123} := \AA^c_0 + \AA^c_1 + \AA^c_1 + \AA^c_2,
$$
where we recall that $ \AA^c_0$ has been defined during the proof of Lemma~\ref{lem:SB1LpLp}.

\smallskip\noindent
{\sl Step 2. Contribution of the term $T_1$. } We prove that 
\beqn\label{eq:estimT1}
T_1 := {1 \over 2} {d \over dt}   \|   f_\BB (t) \|_{L^2(\mu^{1/2})}^2  \le - {1 \over 4} \|   f_\BB (t) \|_{L^2(\nu^{1/2})}^2 +C_1 \,   \|   f_\BB (t) \|_X^2, 
\eeqn
for some positive constant $C_1$ and where $\nu$ is the weight function defined by 
$$
\Blue\nu (x,v) :=  {|v| \over |x|^{1/2}} \, \phi_{1/2}(x). 
$$
In order to prove \eqref{eq:estimT1}, we first observe that 
\bean
T_1 = (\BB \, f_\BB(t),f_\BB(t))_{L^2(\mu^{1/2})},
\eean
and next we compute the RHS by splitting it in several pieces. For any $f \in X$, we have 
\bean
(\BB f, f)_{L^2(\mu^{1/2})} 
&=& (- v \cdot \nabla_x f, f)_{L^2(\mu^{1/2})} +   (\AA^c_{0123} f - Kf, f)_{L^2(\mu^{1/2})}   
\\
&=:& T_{1,1} + T_{1,2}. 
\eean
We compute the first key term. Performing one integration by parts, we find
\bean
T_{1,1} 
&=&\int \bigl( - v \cdot \nabla_x f ) \, f   \Bigl( 1 - {x \over |x|^{1/2}} \cdot {v \over |v|} \Bigr) \, \phi_{1/2}(x) \, dxdv
\\
&=&{1 \over 2} \int f^2 \,  \Bigl\{  \Bigl[  {1 \over 2} \Bigl( {x \over |x|} \cdot {v \over |v|} \Bigr)^2  - 1 \Bigr]    {|v| \over |x|^{1/2}} \phi_{1/2}
+ \Bigl( 1 - {x \over |x|^{1/2}} \cdot {v \over |v|} \Bigr) \, \Bigl( v \cdot \nabla_x \phi_{1/2} \Bigr) \Bigr]\, dxdv
\\
&\le& - {1 \over 4} \int f^2 \,  {|v| \over |x|^{1/2}}  \, \phi_{1/2} \, dxdv + \widetilde T_{1,1} 
\eean
with $  \widetilde T_{1,1}   \lesssim \| f \|_{L^2}^2$. For the remainder term, we also easily have  $ |T_{1,2} |  \lesssim \| f \|_{L^2}^2$ from what 
\eqref{eq:estimT1} immediately follows. 
%
%

\smallskip\noindent
{\sl Step 3. Contribution of the term $T_2$. }  We prove that 
\beqn\label{eq:estimT2}
T_2 :=  {1 \over 2} {d \over dt} \| f_\BB \|^2_X   \le  C_2 \, \|   f_\BB  \|_X^2, 
\eeqn
for some positive constant $C_2$. We have 
\bean
T_2 
&=&  
{1 \over 2} {d \over dt} \| f_\BB  \|^2_{L^2(m)} + {1 \over 2} {d \over dt} \| f_\BB  \|^2_{L^1(m)}
\\
&=&  (\BB \, f_\BB,f_\BB)_{L^2(m)} + \langle \BB \, f_\BB,\hbox{sign} f_\BB \rangle_{L^1(m),L^\infty} \| f_\BB  \|_{L^1(m)}
\\
&=:&  T_{2,1} + T_{2,2},
\eean
with 
\bean
T_{2,1} 
&:=&  (\BB_0 \, f_\BB,f_\BB)_{L^2(m)} + \langle \BB_0 \, f_\BB,\hbox{sign} f_\BB \rangle_{L^1(m),L^\infty} \| f_\BB  \|_{L^1(m)} \le 0,
\eean
because $\BB_0$ is dissipative in both $L^1(m)$ and $L^2(m)$ from Lemma~\ref{lem:SB0Lpm}, 
and with
\bean
T_{2,2} &:=&  (\AA^c_{0123} \, f_\BB,f_\BB)_{L^2(m)} + \langle \AA^c_{0123} \, f_\BB,\hbox{sign} f_\BB \rangle_{L^1(m),L^\infty} \| f_\BB  \|^2_{L^1(m)}
\\
&\le& \| \AA^c_{0123} \, f_\BB \|_{L^2(m)} \, \| f_\BB \|_{L^2(m)}  +\|  \AA^c_{0123} \, f_\BB \|_{L^1(m)}  \| f_\BB  \|^2_{L^1(m)}
\\
&\lesssim&  \| f_\BB \|^2_X,
\eean
\Blue because $| \AA^c_i f | \le A_{1+\chi} |f|$ for any  $i \in \{0, ... ,3 \}$ with $A_{1+\chi} \in \BBB(L^p(m),L^p(m))$ for any $p \in \{1,2\}$. \Black

\smallskip\noindent
{\sl Step 4. Contribution of the term $T_3$. } We prove that for any $\eta_1 \in (0,1)$, we can find $\delta_1,\delta_2,\delta_3 \in (0,1)$ and $R \ge 1$ such that the associated 
operator $\BB$ satisfies 
\beqn\label{eq:estimT3}
T_3 := {1 \over 2} {d \over dt} \int_0^\infty \| S_{\BB_1} (\tau) f_\BB (t) \|^2_X \, d\tau 
   \le  - {3 \over 8}  \, \|   f_\BB (t)  \|_X^2 + {\eta_1 \over 4} \, \|   f_\BB (t)  \|^2_{L^2(\nu^{1/2})}.
\eeqn
 We split the term $T_3$ as 
 \bean
T_3
&=&  
 {1 \over 2} {d \over dt} \int_0^\infty \| S_{\BB_1} (\tau) f_\BB (t) \|^2_{L^1(m)} \, d\tau + {1 \over 2} {d \over dt} \int_0^\infty \| S_{\BB_1} (\tau) f_\BB (t) \|^2_{L^2(m)} \, d\tau
\\
&=:&  T_{3,1} + T_{3,2}. 
\eean
For the first term, we compute 
\bean
T_{3,1}
&=& 
 \int_0^\infty {1 \over 2} {d \over dt} \| S_{\BB_1} (\tau) f_\BB (t) \|^2_{L^1(m)} \, d\tau 
\\ 
&=&
 \int_0^\infty \langle S_{\BB_1} (\tau) \BB f_\BB (t) , \hbox{sign} ( S_{\BB_1} (\tau)  f_\BB (t) ) \rangle_{L^1(m),L^\infty}  \,  \| S_{\BB_1} (\tau) f_\BB (t) \|_{L^1(m)} \, d\tau 
\\ 
&=&
 \int_0^\infty \langle \BB_1 S_{\BB_1} (\tau)  f_\BB (t) , \hbox{sign} ( S_{\BB_1} (\tau)  f_\BB (t) ) \rangle_{L^1(m),L^\infty}  \,  \| S_{\BB_1} (\tau) f_\BB (t) \|_{L^1(m)} \, d\tau 
 \\
&&+
 \int_0^\infty \langle   S_{\BB_1} (\tau)  \AA^c_{123} f_\BB (t) , \hbox{sign} ( S_{\BB_1} (\tau)  f_\BB (t) ) \rangle_{L^1(m),L^\infty}  \,  \| S_{\BB_1} (\tau) f_\BB (t) \|_{L^1(m)} \, d\tau 
\\
&=:&  T_{3,1,1} + T_{3,1,2}. 
\eean
%
On the one hand, we  observe that 
\bean
T_{3,1,1} 
=
  \int_0^\infty {1 \over 2} {d \over d\tau} \| S_{\BB_1} (\tau)   f_\BB (t) \|^2_{L^1(m)} \, d\tau   = - {1 \over 2}   \|   f_\BB (t) \|^2_{L^1(m)},
\eean
where in the last line we have use that $S_{\BB_1} f_0$ has the nice decay estimate \eqref{eq:SB1LpLp} in the space $L^1(m)$ because $f_0 \in X$.  
%
On the other hand, using again the decay estimate \eqref{eq:SB1LpLp},  for any $\eps >0$, we have
\bean
T_{3,1,2}
&=& 
 \int_0^\infty   \| S_{\BB_1} (\tau)  \AA^c_{123} f_\BB (t)  \|_{L^1(m)}  \,  \| S_{\BB_1} (\tau) f_\BB (t) \|_{L^1(m)} \, d\tau 
\\
&\le&
 {1 \over 2\eps} \int_0^\infty   \| S_{\BB_1} (\tau)  \AA^c_{123} f_\BB (t)  \|_{L^1(m)}^2 \, d\tau 
 +  {\eps \over 2} \int_0^\infty   \,  \| S_{\BB_1} (\tau) f_\BB (t) \|_{L^1(m)} ^2 \, d\tau 
\\
&\lesssim&
{1 \over \eps}    \|  \AA^c_{123} f_\BB (t)  \|_X^2 
 + \eps  \,  \|  f_\BB (t) \|_X^2.
\eean
We may treat the second term in a similar way, using in particular the fact that $S_{\BB_1} f_0$ has the nice decay estimate \eqref{eq:SB1LpLp} in the space $L^2(m)$ because $f_0 \in X$. Gathering the two resulting estimates and taking $\eps > 0$ small enough, we get  
\beqn\label{eq:estimT3bis}
T_3
\le - {7 \over 16}   \|   f_\BB (t) \|^2_X + \widetilde {T_3}, \quad 
 \widetilde {T_3} \lesssim  \| \AA^c_{123}  f_\BB (t) \|_X^2. 
\eeqn
In order to conclude, we compute the contributions $ \| \AA^c_{i}  f_\BB (t) \|_X^2$ for any $i \in \{1,2,3\}$. 

\smallskip
On the one hand, using the Cauchy-Schwarz inequality, for any $p \in \{1,2\}$, $R > 2$ and $\delta_1 \in (0,1)$, we have 
\bean
\| \AA^c _1 f \|_{L^p}^p(m) 
&=& \int \Bigl|  \phi_R(x)  \int_\VV K' f' \psi^c_{\delta_1}(v') dv' \Bigr|^p  m^p \, dxdv
\\
&\le& 2^p \, m(2R)^p \int   \Bigl[ \int_\VV  | f' |^2  \, dv' \Bigr|^{p/2}  \Bigl[ \int_\VV  \bigl( {\bf 1}_{|v'| \le 2{\delta_1}}  +  {\bf 1}_{V_0 - |v'| \le 2{\delta_1}} \bigr)\, dv' \Bigr]^{p/2}   \, dxdv
\\
&\lesssim&  m(2R)^p \, \delta_1^{p/2} \, \| f \|_{L^2}^{p/2}.
\eean
Similarly, when furthermore $\delta_2 \in (0,1/4)$, we have 
\bean
\| \AA^c _2 f \|_{L^p}^p(m) 
&=& \int \Bigl|  \phi_{\delta_2}(x)  \int_\VV K' f' \psi_{\delta_1}(v') dv' \Bigr|^p  m^p \, dxdv
\\
&\le& 2^p \, m(2)^p \Bigl[  \int_{B(0,1) \times \VV}   {\bf 1}_{|x| \le {2\delta_2}}      | f' |^2 \,  {\bf 1}_{|v'| \ge {\delta_1}}  \, dv' dx  \Bigr]^{p/2}  
\\
&\lesssim& \Bigl[  {\delta_2^{1/2} \over \delta_1} \Bigr]^{p/2} \, \| f \|^p_{L^2(\Blue \nu^{1/2})}. 
\eean
Finally and similarly again, when  furthermore $\delta_3 \in (0,1/2)$, observing that 
$$
0 \le K^c_{\delta_3}(x,v') = \chi \, ( \zeta - \zeta_{\delta_3}) (x \cdot v') \le  \chi \, {\bf 1}_{|x \cdot v'| \le \delta_3},
$$
we have 
\bean
\| \AA^c _3 f \|_{L^p}^p(m) 
&=& \int \Bigl| \phi_{\delta_2,R}(x)  \int K_{\delta_3}^c(x \cdot v') f' \psi_{\delta_1}(v') \, dv' \Bigr|^p  m^p \, dxdv
\\
&\le& m(2R)^p \, \chi^p \int   \Bigl[ \int_\VV  | f' |^2  \, dv' \Bigr|^{p/2}  \Bigl[ \int_\VV  {\bf 1}_{|x \cdot v'| \le 2{\delta_3}}  \, dv' \Bigr]^{p/2}   \, {\bf 1}_{|x| \ge \delta_2} \, dxdv
\\
&\lesssim&  m(2R)^p  \,  \Bigl[   \int   f^2  \, dv dx  \Bigr|^{p/2} \Bigl[ \hbox{meas} \bigl\{ v \in \VV; \, |v_1| \le \delta_3/\delta_2 \bigr\} \Bigr]^{p/2} 
\\
&\lesssim&  m(2R)^p  \, {\delta_3^{p/2} \over \delta^{p/2}_2} \, \| f \|_{L^2}^{p} .
\eean
All these estimates together, we get 
\beqn\label{eq:estimT3ter}
 \| \AA^c_{123}  f_\BB (t) \|_X^2
 \lesssim   m(2R)^p \, \delta_1 \, \|  f \|^2_X +  {\delta_2^{1/2} \over \delta_1 }   \, \| f \|^2_{L^2(\Blue \nu^{1/2})} + m(2R)^p  \, {\delta_3  \over \delta_2} \, \| f \|^2_X  .
 \eeqn
 We thus obtain \eqref{eq:estimT3} by just gathering \eqref{eq:estimT3bis} and \eqref{eq:estimT3ter} and by choosing $ \delta_1, \delta_2, \delta_3 > 0$ adequately.

\smallskip\noindent
{\sl Step 5.  Conclusion.}
%
%
\Blue From  estimates \eqref{eq:equivNX}, \eqref{eq:defT=T123}, \eqref{eq:estimT1},  \eqref{eq:estimT2} and \eqref{eq:estimT3}, we have 
\bean
T  
&\le&  \eta_1 \, C_1 \, \| f_\BB \|_{X}^2 +  \eta_2 \, C_2 \,  \| f_\BB \|_X^2  - {3 \over 8} \| f_\BB \|^2_X
\\
 \\
&\le& -  {1 \over 4} \| f_\BB \|^2_X \le a^* N(f_\BB)^2,
\eean 
by choosing $\eta_1,\eta_2 > 0$  small enough. We have proved that $\BB-a^*$ is dissipative \Black for the norm $N(\cdot)$ and thus \eqref{eq:SBX} follows. \qed

 \subsection{Some regularity associated to  $\AA S_{\BB}$ }
\label{subseq:REstimASB}
In this section we show that the family of operators $\AA S_\BB$ satisfies a regularity and growth estimate
that we express in terms of the abstract Sobolev space $X^{1/2}_\BB$ defined as the usual $1/2$ interpolated space between $X$ and
the domain
$$
X^1_\BB = D(\BB) := \{ f \in X; \BB f \in X \}
$$
endowed with the graph norm.

 \begin{prop}\label{prop:Reg3}  For the same constant  $a^* < 0$ as defined in Lemma~\ref{lem:SB1L1m}, for any $a > a^*$ there exits $C_a \in (0,\infty)$ such that the family of operators $\AA S_\BB$ satisfies  
\beqn\label{eq:Reg3SB}
\int_0^\infty \| \AA S_\BB  (t)  \, f \|^2_{ Y}  \, e^{- 2a t} \, dt \le C_a \,  \| f \|^2_X, \quad \forall \, f \in X,
\eeqn
with 
$$
Y :=  \{ f \in L^2(\R^d \times \VV); \,\, \hbox{supp} \, f \subset B(0,R) \times \VV, \,\, f \in H^{1/2} \}.
$$
\end{prop}

The proof is mainly a consequence of Bouchut-Desvillettes' version \cite[Theorem 2.1]{BouchutD} (see also \cite{DesvM} for a related discrete version) of the classical averaging Lemma initiated in the famous articles of Golse et al. \cite{GolsePS,GolseLPS}. We give in step 1 below a simpler, more accurate and more adapted version of  \cite[Theorem 2.1]{BouchutD} for which we sketch the proof for  the sake of completeness. During the proof, we will use the following classical trace result. 
\begin{lem}\label{lem:ch?:TraceResult}
There exists a constant $C_d \in(0,\infty)$ such that for any $\phi \in H^{d/2} (\R^d)$ and any $u \in \R^d$, $|u| =1$, the \Blue real function $\phi_u$, defined by $\phi_u(s) := \phi (s u)$
\Blue for any $s \in \R$, \Black  satisfies
$$
\| \phi_u \|_{L^2(\R)}\le C_d \| \phi  \|_{H^{d/2}(\R^d)}  = C_d \left( \int_{\R^d} |{ \check F \phi }|^2 (w) \, \langle w \rangle^d \, dw \right)^{1/2},
$$
where  $\check F$ stands for the (inverse) Fourier transform operator. 
\end{lem}

\noindent
{\sl Proof of Proposition~\ref{prop:Reg3}. } We split the proof into two steps. 

\smallskip\noindent
{\sl Step 1. } We consider the damped free transport equation 
\beqn\label{eq:FreeTranspEq}
\partial_t f = \TT f := -  v \cdot \nabla_x f - f, \quad f_{|t=0} = f_0,
\eeqn
and we denote by $S_\TT(t)$ the associated semigroup defined through the characteristics formula
\beqn\label{eq:FreeTranspCF}
[S_\TT(t) f_0] (x,v] := f(t,x,v) = f_0 (x-vt,v) \, e^{-t}. 
\eeqn
We claim that for any $\varphi \in L^2(\VV)$, there holds
\beqn\label{eq:ASSL2H12}
\int_0^\infty \| A_\varphi S_\TT (t) \varphi \|_{H^{1/2}_{x}}^2 \, e^{2t} \,  dt  \lesssim \| \varphi \|_{L^2(\VV)}.  
\eeqn
For a given function $h$ which depends on the $x$ variable or on the $(x,v)$ variable, 
 we denote by $\hat h$ its Fourier transform on the $x$ variable and by $\FF h$ its Fourier transform on both variables $x$ and $v$.  
We fix $f_0 \in L^2(\R^d \times \VV)$ and $\varphi \in L^\infty(\R^d)$, we denote by $f$ the solution to the free transport equation 
\eqref{eq:FreeTranspEq} and by $\rho$ the average function
$$
\rho(t,x) := \int_{\R^d} f(t,x,v)  \, \varphi (v) \, dv = [A_\varphi S_\TT(t) f_0 ] (x).
$$
In Fourier variables, the free transport equation \eqref{eq:FreeTranspEq} writes 
$$
\partial_t \hat f + i v \cdot \xi \hat  f - \hat f= 0, \quad \hat f_{|t=0} = \hat f_0,
$$
so that 
$$
\hat f (t,\xi,v) = e^{i v \cdot \xi \, t - t} \hat f_0 (\xi, v)
$$
and
\bean
\hat \rho (t,\xi) 
= \int_{\R^d} e^{i v \cdot \xi \, t -t} \hat f_0 (\xi, v) \, \varphi (v) \, dv
=  \FF ( f_0  \, \varphi ) (\xi, t \xi) \, e^{-t}.
\eean
We deduce 
$$
\int_0^\infty |\hat \rho (t,\xi) |^2 \, e^{2t} \, dt \le \int_\R | \FF ( f_0  \, \varphi ) (\xi, t \xi) |^2 \, dt.  
$$

Performing one change of variable, introducing the notation $\sigma_\xi = \xi/|\xi|$ and using Lemma~\ref{lem:ch?:TraceResult}, we deduce 
\bean
\int_\R | \FF ( f_0  \, \varphi ) (\xi, t \xi) |^2 \, dt 
&=& {1 \over |\xi|} \int_\R | \FF ( f_0  \, \varphi ) (\xi, s \, \sigma_\xi) |^2 \, ds 
\\
&\lesssim& {1 \over |\xi|} \int_{\R^d} | (\hat f_0  \, \varphi ) (\xi, w) |^2 \langle w \rangle^d \, dw.
\eean
Thanks to Plancherel identity, we then obtain 
\bean
\int_0^\infty \!\!  \int_{\R^d} |\xi| \, |\hat \rho (t,\xi) |^2 \, d\xi \, e^{2t}  dt 
&\lesssim& \int_{\R^d} \!\!  \int_{\R^d} |  (f_0  \, \varphi ) (x, w) |^2 \langle w \rangle^d \, dw  dx
= \| \varphi \|_{L^2_{d/2}}^2 \, \| f_0 \|^2_{L^2_{xv}}, 
\eean
which ends the proof \eqref{eq:ASSL2H12}.

\smallskip\noindent
{\sl Step 2. We show a similar estimate on $\AA S_\TT(t)$. } Using that $ K_{R,\delta_i} \in C^\infty_c(\R^d \times \R^d)$, supp$\, K_{R,\delta_i} \subset B(0,2R) \cap B(0,V'_0)$, $V'_0 \in (0,V_0)$, we may expand it as a Fourier series 
$$
 K_{R,\delta_i} (x,v) = \sum_{k,\ell \in \Z^{d}} a_{k,\ell} \, e^{i \, x \cdot k} \, e^{i \, v \cdot \ell} \, \vartheta (v), \quad \forall \, (x,v) \in \QQ,
$$
$\QQ := \{ x \in \R^d, \, v \in \R^d; \,\, |x_i| \le 2R, \, |v_i | \le V_0, \,\, \forall \, i = 1, ..., d \}$, 
for a truncation function $\vartheta\in C^\infty(\R^d)$, supp$\, \vartheta \subset B(0,V_0)$, $\vartheta \equiv 1$ on $B(0,V'_0)$ and with fast decaying Fourier coefficients   
$$
|a_{k,\ell} | \lesssim \langle k \rangle^{-2d-4} \,  \langle  \ell \rangle^{-2d-2}.
$$
From the definition of $\AA$ and denoting $f (t) = S_\TT (t) \, f_0$ for some $f_0 \in L^2(\R^d \times \VV)$, we may then write 
$$
(\AA S_\TT (t) f_0 )(x) = \sum_{k,\ell \in \Z^{d}} a_{k,\ell} \, e^{i \, x \cdot k} \, \rho_\ell (t,x),
\quad \rho_\ell(t,x) := \int_\VV  f(t,x,v) \, e^{i \, v \cdot \ell} \, \vartheta (v) \, dv.
$$
On the one hand, from Step 1, we have 
\beqn\label{eq:Reg3Step2Ineq1}
\sup_{\ell \in \Z^d} \int_0^\infty \| \rho_\ell(t, \cdot) \|^2_{H^{1/2}} \, e^{2t} \, dt \lesssim \| f_0 \|_{L^2}^2.
\eeqn
On the other hand, we denote $e_k(x) := e^{i \, x \cdot k}$ and we define the mapping  
$$
 U(\rho_\ell) :=  \sum_{k,\ell \in \Z^{d}} a_{k,\ell} \, e_k \, \rho_\ell.
$$
From Cauchy-Schwarz inequality and Fubini Theorem, we have 
\bean
&&\int_0^\infty \| U(\rho_\ell)(t,\cdot) \|_{L^2(B_{2R})}^2\, e^{2t}  \, dt \le
\\
&&\quad\le \int_0^\infty \int_{B_{2R}} \Bigl( \sum_{k,\ell} |a_{k,\ell}|^2 \,  \langle k \rangle^{d+1} \,  \langle  \ell \rangle^{d+1} \Bigr) 
\Bigl( \sum_{k,\ell} |\rho_\ell |^2 \,  \langle k \rangle^{-d-1} \,  \langle  \ell \rangle^{-d-1} \Bigr) \, e^{2t}  \, dxdt
\\
&&\quad\lesssim \sum_{k,\ell}  \langle k \rangle^{-d-1} \,  \langle  \ell \rangle^{-d-1} \int_0^\infty \int_{B_{2R}}   
 |\rho_\ell |^2 \, e^{2t}    \, dtdx
\\
&&\quad\lesssim \sup_{\ell \in \Z^d} \int_0^\infty  \| \rho_\ell (t,\cdot) \|_{L^2(B_R)}^2 \, e^{2t}  \, dt. 
\eean
Using furthermore that 
$$
 \nabla_x U(\rho_\ell) =   \sum_{k,\ell \in \Z^{d}} a_{k,\ell} \, (ik) \, e_k \, \rho_\ell +  \sum_{k,\ell \in \Z^{d}} a_{k,\ell} \, e_k \, \nabla_x\rho_\ell ,
$$
we find similarly
\bean
\int_0^\infty \| \nabla_x U(\rho_\ell)(t,\cdot) \|_{L^2(B_{2R})}^2 \, e^{2t} \, dt\lesssim \sup_{\ell \in \Z^d} \int_0^\infty  \| \rho_\ell (t,\cdot) \|_{H^1(B_R)}^2\, e^{2t}  \, dt.
\eean
Observing that 
$$
\{ g \in L^2(\R^d \times \VV); \, \hbox{supp} \, g \subset B(0,2R) \times \VV, \, \nabla_x g \in L^2 \} \subset X^1_\BB,
$$
both estimates together and an interpolation argument yield
\beqn\label{eq:Reg3Step2Ineq2}
\int_0^\infty \| U(\rho_\ell)(t,\cdot) \|_{X^{1/2}_\BB}^2 \, e^{2t}  \, dt\lesssim \sup_{\ell \in \Z^d} \int_0^\infty  \| \rho_\ell (t,\cdot) \|_{H^{1/2}(B_R)}^2 \, e^{2t} \, dt.
\eeqn
Gathering estimates \eqref{eq:Reg3Step2Ineq1} and \eqref{eq:Reg3Step2Ineq2}, we have established 
\beqn\label{eq:Reg3Step2Ineq3}
\int_0^\infty \| \AA S_\TT (t) f_0 \|_{X^{1/2}_\BB}^2 \, e^{2t}  \, dt\lesssim   \| f_0 \|_{L^2}^2.
\eeqn

\smallskip\noindent
{\sl Step 3. Conclusion. }  We split $\BB$ as $\BB = \TT+ \CC$. The Duhamel formula writes
$$
S_\BB = S_\TT + S_\TT * \CC S_\BB,
$$
from which we deduce 
$$
\AA S_\BB = \AA S_\TT + \AA S_\TT * \CC S_\BB. 
$$
We just have to bound the last term in order to establish  \eqref{eq:Reg3SB}. For that purpose, we fix $f \in X$, $a > \alpha > a^*$, and we compute
\bean
\int_0^\infty \| \AA S_\TT * \CC S_\BB (t) f \|_Y^2 \, e^{- 2 a t } \, dt
&\le&
\int_0^\infty\int_0^t \| \AA S_\TT (t-s) \CC S_\BB (s) f \|_Y^2 \,  ds \, t\, e^{- 2 a t }\, dt
\\
&\le&
\int_0^\infty\int_0^\infty \| \AA S_\TT (\tau) \CC S_\BB (s) f \|_Y^2 \,   e^{- 2 \alpha \tau }\, d\tau e^{- 2 \alpha s }ds
\\
&\lesssim&
\int_0^\infty  \|   \CC S_\BB (s) f \|_X^2  e^{- 2 \alpha s }ds
\lesssim \| f \|_X^2,
\eean
where we have used the Cauchy-Schwarz  inequality, estimates  \eqref{eq:Reg3Step2Ineq3} and \eqref{eq:SBX}. \qed

%
\subsection{A first asymptotic stability estimate in $X$}
 \label{subsec:ASEinX}
 In order to apply the semigroup version \cite[Theorem 5.3]{MiSch2016} and \cite{MSerratum*} of the Krein-Rutman theorem, 
 we list below some properties satisfied by the operators $\LL$, $\AA$ and $\BB$. 
 
\smallskip\noindent
{\sl Fact 1. }  There exists $a^* < 0$ such that for any $a>a^*$ and $\ell\in\N$, 
 the following growth estimate holds
$$
t \mapsto    \| (\AA S_\BB )^{(\ast\ell)}(t)\|_{\BBB(X)} \, e^{-a t} \in L^\infty(\R_+). 
$$
That is an immediate consequence of Proposition~\ref{prop:SBLpLp} and $\AA\in\BBB(X)$. 

\smallskip\noindent
{\sl Fact 2. }  We define the resolvent operator 
$$
R_\BB(z) := (\BB-z)^{-1} = - \int_0^\infty S_\BB(t) \, e^{-zt} \, dt
$$
for $z \in \Delta_a := \{ \zeta \in \C; \, \Re e \, \zeta  > a \}$ and $a$ large enough. For the same value $a^* < 0$ as above, there exists $Y \subset X^s_\LL$, $s \in (0,1/2)$, with compact embedding such that  for any $a>a^*$  the following  estimate holds
$$
 \|\AA R_\BB (z)\|_{\BBB(X,Y)} \le C_a, \quad\forall \, z \in \Delta_a .
$$
That is an immediate consequence of  Proposition~\ref{prop:Reg3}, which readily implies
\bean
\| \AA R_\BB(z) f \|^2_Y 
&\le& \int_0^\infty \| \AA S_\BB(t) f\|^2_Y \langle t \rangle^2 \, e^{-2at} \, dt \int_0^\infty \langle t \rangle^{-2} \, dt
\\
&\lesssim& \| f\|^2_X, \quad \forall \, f \in X, 
\eean
together with the fact that 
$$
 \{ f \in L^2(\R^d \times \VV); \,\, \hbox{supp} \, f \subset B(0,R) \times \VV, \,\, f \in H^1 \} \subset  X^1_\LL
$$
and an interpolation argument. 
 
\smallskip\noindent
\Blue {\sl Fact 3. } The semigroup $S_\LL$ is positive, the operator $-\LL$ satisfies the strong maximum principle as stated in Lemma~\ref{lem:-Lspos}
and $\LL$ satisfies Kato's inequality \eqref{eq:KatoIneq}. \Black
  
\smallskip\noindent
{\sl Fact 4. } The mass conservation property writes $\LL^* 1 = 0$, so that $0 > a^*$ and $0$  is an eigenvalue for the dual problem
associated to a positive dual eigenfunction. 

\smallskip
Gathering these above facts, we may then apply \cite[Theorem 5.3]{MiSch2016}, or more exactly we may repeat the proof of \cite[Theorem 5.3]{MiSch2016} with minor and straightforward adaptations (we refer to \cite{Mbook*} where these slight modifications are performed), in order to obtain that $0$ is a (algebraically) simple eigenvalue, that there exists a spectral gap between this largest eigenvalue $0$ and the remainder part of the spectrum and that a  quantitative (partial but principal) spectral mapping theorem holds true. More precisely, we have the following asymptotic estimate: there exists $\alpha \in (a_*,0)$ such that  
\beqn\label{eq:theo2X}
\| S_\LL(t) \Pi^\perp f_0 \|_{X} \lesssim  e^{at} \,  \|   f_0  \|_{X},
\quad \forall \, f_0 \in X, \,\,  \forall \, a > \alpha, \,\,  \forall \, t \ge 0,
\eeqn
where we have set $\Pi^\perp := I - \Pi$ and $\Pi f_0 := \langle \!\langle f_0 \rangle\!\rangle G$. 

%
\subsection{Asymptotic stability estimate in weighted $L^1$ spaces}
\label{subseq:ASEinL1} 

We first consider the exponential weight $m(x) := \exp(\gamma \langle x \rangle)$ with $\gamma \in (0,\gamma^*)$ and $\gamma^* > 0$ identified in Lemma~\ref{lem:estimL1m}. 
Iterating the Duhamel formula, we may write  
$$
S_\LL \Pi^\perp = \Pi^\perp \{ \SS_{\BB_1} +  ... + \SS_{\BB_1} * (\AA_1  S_{\BB_1})^{N-1} \} + (\SS_{\LL} \Pi^\perp )* (\AA_1  S_{\BB_1})^{N},
$$
with $N=d+2$. From Lemma~\ref{lem:SB1L1m}  and \eqref{eq:ASB1L1toLinfty} we have $\SS_{\BB_1} * (\AA_1  S_{\BB_1})^{\ell} : L^1(m) \to L^1(m)$ with rate $e^{at}$ for any $\ell \in \{ 0, ... , N-1\}$ and  $(\AA_1  S_{\BB_1})^{N} : L^1 (m) \to X$ with rate $e^{at}$. Using that 
$\SS_{\LL} \Pi^\perp : X \to X \subset L^1(m)$ with rate $e^{at}$ from \eqref{eq:theo2X} and gathering all the preceding decay estimates, we conclude that \eqref{eq:theo2} holds in $L^1(m)$.


\medskip
We next consider the polynomial weight   $m(x) := \langle x \rangle^k$ with $k \in (0,\infty)$.  
We begin with a decay estimate on the semigroup $S_{\BB_1}$.

\begin{lem}\label{lem:SB1L1k}
For any  $k>\ell > 0$, 
the semigroup
$S_{\BB_1}$ satisfies the following growth estimate
$$
\|S_{\BB_1}\|_{L^1_k\to L^1_\ell}\lesssim \langle t\rangle^{- (k-\ell)},
\quad\forall t\ge0.
$$
\end{lem}

\noindent
{\sl Proof of Lemma~\ref{lem:SB1L1k}. } Recalling that the dual operator $\LL^*$ has been defined in the proof of Lemma~\ref{lem:estimL1m}, for any $q > 0$, we compute   
\bean
\LL^*\langle\gamma x \rangle^q = q\gamma (v \cdot x) \langle \gamma x \rangle^{q-2},
\eean
\bean
\LL^*(v \cdot x) \langle\gamma x \rangle^{q-2} 
&=& v \cdot \nabla_x [ (v \cdot x) \langle\gamma x \rangle^{q-2}] - q \, (v \cdot x) \langle\gamma x \rangle^{q-2}
\\
&=& \Bigl( |v|^2 -(v \cdot x)-\chi |v \cdot
x|\Big)\langle\gamma x \rangle^{q-2}
+(q-2)\gamma(v \cdot x)^2 \langle\gamma x \rangle^{q-4}.
\eean
We then compute 
\bean
\LL^*|v \cdot x|\langle\gamma x\rangle^{q-2} &=&  
v \cdot \nabla_x [|v \cdot x|\langle\gamma x\rangle^{q-2}] +
( 1 + \chi \zeta) \Big(V_1 |x|\langle\gamma x\rangle^{q-2} - |v\cdot
x|\langle\gamma x\rangle^{q-2}\Big)
\\
&=& \Bigl( |v|^2 \, {v \cdot x \over |v \cdot x|} +( 1 + \chi
\zeta) V_1 |x| -  ( 1 + \chi \zeta)  |v\cdot x| \Bigr) \langle\gamma x\rangle^{q-2}\\
&& +(q-2)\gamma|v\cdot x|(v\cdot x)\langle\gamma x\rangle^{q-4} ,
\eean
where we recall that $V_1$ has been defined in \eqref{eq:defV1}. 
We consider $\beta,\gamma \in (0,1)$ to be fixed later  such that the weight function
$$
\widetilde m_q :=\langle\gamma x\rangle^q+q\gamma(v\cdot x)\langle\gamma x\rangle^{q-2}
-q\beta|v\cdot x|\langle\gamma x\rangle^{q-2}
$$
satisfies 
$$
(1-\delta)\langle\gamma x \rangle^q \, \le  \widetilde m_q \le (1+\delta) \,\langle\gamma x \rangle^q ,  
$$
for some $\delta \in (0,1)$. Gathering the previous estimates, there holds
\bean
\BB_1^*\widetilde m_q &=& \LL^*\widetilde m_q-\AA_1^*\widetilde m_q= 
\LL^*\widetilde m_q -(1+\chi\zeta)\phi_R\int_\VV \widetilde m_q\mathrm{d}v
\\
&=& q\gamma|v|^2\langle\gamma
x\rangle^{q-2}+q(q-2)\gamma (v\cdot x)^2\langle\gamma
x\rangle^{q-4}-q\gamma\chi|v\cdot x|\langle\gamma x\rangle^{q-2}
\\
&& -q\beta |v|^2\frac{v\cdot x}{|v\cdot x|}\langle\gamma
x\rangle^{q-2}-q(q-2)\beta\gamma|v\cdot x|(v\cdot
x)\langle\gamma x\rangle^{q-4}
\\
&& -q\beta (1+\chi\zeta)V_1|x|\langle\gamma
x\rangle^{q-2}\phi^c_R+q\beta (1+\chi\zeta)|v\cdot
x|\langle\gamma x\rangle^{q-2}
\\
&& -(1+\chi\zeta)\langle\gamma x\rangle^q\phi_R,
\eean
and then 
\bean
\BB_1^*\widetilde m_q 
&\le& q\gamma V_0^2\langle\gamma x\rangle^{q-2}+q|q-2|\gamma 
V_0^2|x|^2\langle\gamma x\rangle^{q-4}-q\gamma\chi |v\cdot
x|\langle\gamma x\rangle^{q-2}
\\
&& +q\beta V_0^2 \langle\gamma x\rangle^{q-2}+q|q-2|\beta\gamma 
V_0^2 |x|^2\langle\gamma x\rangle^{q-2}
\\
&& -q\beta(1-\chi)V_1|x|\langle\gamma
x\rangle^{q-2}+q\beta(1+\chi) |v\cdot x|\langle\gamma
x\rangle^{q-2}
\\
&& -(1-\chi) \langle\gamma
x\rangle^{q-1}\phi_R+q\beta(1-\chi)V_1 \langle\gamma x\rangle^{q-1}\phi_R
\\
&\le& \Big(q V_0^2(\gamma+\beta)+q\beta(1-\chi)V_1\Big)\langle\gamma
x\rangle^{q-2}+q|q-2|V_0^2\gamma(1+\beta)|x|^2\langle\gamma 
x\rangle^{q-4}
\\
&&-(1-\chi)(1-q\beta V_1)\langle\gamma
x\rangle^q\phi_R-q\beta(1-\chi)V_1\langle\gamma x\rangle^{q-1}
\\
&\le& \Big(\frac{C_1}{\langle\gamma x\rangle}-C_2\phi_R
-q\beta(1-\chi)V_1\Big)\langle\gamma x\rangle^{q-1}.
\eean
Choosing $\beta(1+\chi)=\gamma\chi$ with $\gamma>0$ small enough and $R\ge1$ large enough, and observing that  $C_1=O(\gamma)$, $C_2 \ge (1-\chi)/2$ as $\gamma \to 0$,  
we deduce 
$$
\BB_1^*\widetilde m_q 
\le - {q\beta(1-\chi)V_1 \over 2} \langle\gamma x\rangle^{q-1} 
\lesssim  - \langle  x\rangle^{q-1}.
$$
We denote $f_{\BB_1}(t) :=S_{\BB_1}(t) f_0$ for some $0 \le f_0 \in L^1_k$ and then $\widetilde M_q = \langle\!\langle f_{\BB_1} \, \widetilde m_q \rangle\!\rangle$, 
$M_q = \langle\!\langle f_{\BB_1} \langle x \rangle^q \rangle\!\rangle$, so that 
\beqn\label{eq:MqtildeMq}
\widetilde M_q \lesssim M_q \lesssim \widetilde M_q.
\eeqn From the above inequality, we get
\beqn\label{eq:ddtMtilde}
{\mathrm{d}\over\mathrm{d}t} \widetilde M_q =\int f \, (\BB_1^*\widetilde m_q)
 \lesssim  - M_{q-1},
\eeqn
and in particular 
\beqn\label{eq:Mtildekt0}
\widetilde M_k(t) \le \widetilde M_k(0), \quad \forall \, t \ge 0.
\eeqn
A classical interpolation inequality together with \eqref{eq:MqtildeMq} and \eqref{eq:Mtildekt0} give 
$$
M_\ell(t) \le M_{\ell-1}(t)^{\theta} M_{k}(t)^{1-\theta} 
\lesssim M_{\ell-1}(t)^{\theta} M_{k}(0)^{1-\theta},
$$
with $\theta \in (0,1)$ such that $\ell = \theta (\ell-1) + (1-\theta) k$. Coming back to \eqref{eq:ddtMtilde}, we get 
$$
\frac{\mathrm{d}}{\mathrm{d}t} \widetilde M_\ell  \lesssim - M_k(0)^{- 1/\alpha} \widetilde M_\ell ^{1 + 1/\alpha}
$$
where
$$
\alpha := {1 \over (1/\theta)-1} =  k-\ell.
$$
Integrating the above differential inequality, we obtain 
$$
M_\ell(t) \lesssim {M_{k}(0) \over t^\alpha}, \quad \forall \, t > 0,
$$
and we conclude gathering that last inequality with \eqref{eq:Mtildekt0}. 
\qed

\medskip
In order to establish the asymptotic stability in $L^1(m)$ for a polynomial weight $m$,
we write 
$$
S_\LL \Pi^\perp = \Pi^\perp \SS_{\BB_1} + (\SS_{\LL} \Pi^\perp )* (\AA_1  S_{\BB_1}). 
$$
Introducing the exponential weight  $m_0 := e^{\langle x \rangle}$, we observe that $\Pi^\perp \SS_{\BB_1} : L^1(m) \to L^1$ and
$\AA_1  S_{\BB_1} L^1(m) \to L^1(m_0)$, with rate $\langle t \rangle^{-\ell}$ for any $\ell \in (0,k)$ from Lemma~\ref{lem:SB1L1k}. 
Because we have already established that $\SS_{\LL} \Pi^\perp : L^1(m_0) \to L^1$ with rate $e^{at}$ for any $a \in (a^*,0)$, 
we immediately conclude that \eqref{eq:theo2} holds in $L^1(m)$.

\bigskip
\bibliographystyle{acm}

\def\cprime{$'$}


\end{document}